# A VMS-FEM for the stress-history-dependent materials (fluid or solid) interacting with the rigid body structure

# Part A:

# formulation, numerical verification and application in the pipe-soil-water interaction analysis

**First author:**

Tianyu Li

**Affiliation:**

University of Oxford

# A VMS-FEM for the stress-history-dependent materials (fluid or solid) interacting with the rigid body structure

# Part A:

# formulation, numerical verification and application in the pipe-soil-water interaction analysis

**Abstract:**

The problem of multiphase materials (fluid or solid) interacting with the rigid body structure is studied by proposing a novel VMS-FEM (variational multi-scale finite element method) in the Eulerian framework using the fixed mesh. The incompressible N-S equation with high Reynolds number is stabilized using the idea of the VMS stabilization. To model the multiphase materials, a modified level set method is used to track the development of interfaces between different phases of materials. The interaction between the materials and the rigid body structure is realized using the Lagrangian multiplier method. The mesh of the rigid body structure domain is incompatible with the Eulerian fixed mesh for material domains. A cutting finite element scheme is used to satisfy the coupling condition between the rigid body structure boundary and the deformeable multiphase materials.

A novel algorithm is proposed to track the stress history of the material point in the Eulerian fixed-mesh framework such that the stress-history-dependent fluid or solid materials (e.g. as soil) are also analysable using fixed-mesh. The non-local stress theory is used. To consider the stress-history-dependent material in the Eulerian fixed mesh, in the step of calculating the non-local relative-distance-based averaging treatment, the field point location is based on the updated material point position, rather than the spatial field point location. Thus, the convective effect of the previous stress history and the calculation of the non-local stress fields are realized in one shot. Thus, the proposed novel VMS-FEM can solve almost any history-dependent materials (fluid or solid) in the Eulerian framework (fixed mesh), where the distorted mesh problem is not existent. Moreover, the complicated contact phenomena between different phases of materials can be simulated automatically without using any dedicated contact algorithm as required in the Lagrangian framework using the material-attached mesh. The elasplastic constitutive model with the strain-rate-dependent hardening effect and the plastic-strain-dependent softening effect is derived in the Eulerian framework to model the soil. The interaction problem between the rigid pipe and the soil having very large deformation and strain is studied. The proposed VMS-FEM shows very good accuracy in the benchmark problems. The pipe-soil-water interaction analysis is presented.

# 1. Introduction

A novel VMS-FEM is proposed for the interaction problem between the rigid structure and the multiphase materials (fluid or solid) with history dependent stress (e.g. the elasplastic constitutive model with the plastic-strain-rate-dependent hardening effect and the plastic-strain-dependent softening effect) in the Eulerian framework using the fixed-mesh. Different from the Lagrangian framework where the mesh deforms with the material, in the Eulerian framework, the mesh is fixed into the space and the materials flow through the fixed mesh. The Lagrangian framework and Eulerian frameworks both have their advantages and disadvantages. The Lagrangian framework can track the stress history conveniently because the mesh is attached with the material. But the distorted mesh in the large deformation analysis and the complicated contact phenomena between different phases of materials are two main problems. The self-adaptive mesh and the dedicated contact algorithm must be used; the simulation will be very time consuming for the nonlinear problems, which also need the mapping calculation of stress histories between different time step domains.

The Eulerian framework can study the contact phenomena between different phases of materials automatically. The mesh is fixed in space; the materials flow through the fixed mesh. Thus, the contact condition between different phases of materials can be satisfied automatically. The distorted mesh problem is also not existent in the Eulerian framework because the mesh is fixed.

However, the stress-history-dependent materials cannot be analysed directly in Eulerian framework. To track the stress history in the Eulerian framework is not trivial because the materials flow through the fixed mesh. In this paper, a novel VMS-FEM is proposed for the N-S equation in the Eulerian framework such that the stress-history-dependent materials (solid or fluid) can also be studied in the way of computational fluid mechanics. The contact phenomena between different phases of materials is perfectly satisfied automatically. The distorted mesh problem is also not existent. There is no need to perform the self-adaptive mesh. Thus, the advantages of the Eulerian and Lagrangian frameworks are both included in this method.

The N-S equation is the governing equation for the fluid mechanics. But actually, it is also able to study the solid mechanics. There is no need to distinguish the solid mechanics and the fluid mechanics in numerical methods. As reported in (Liu and Marsden, 2018), the solid mechanics and fluid mechanics can be studied in a unified framework. The solid mechanics is different from the fluid mechanics because of its stress-history-dependency, meaning that the current stress increment is not only dependent on its current deformation status, but also dependent on the history of the accumulated plastic strain, the previous stress history and the potential internal variable history. To consider solid mechanics using N-S equation, one extra stress term is added to represent the history-dependent stress. The corresponding VMS-FEM formulation is derived in this paper.

The non-local stress theory is a relative-distance-weighted averaging algorithm that can clean the noise of the original stress field. Typically, the source point's stress is calculated as a relative-distance-weighted summation of the stress values of all the nearby field points. The relative distance between the source point and the field point determines the weight function uniquely.

The calculation of the non-local stress is trivial in the Lagrangian framework, where the mesh and the material are always bonded. However, in the Eulerian framework, the stress of a given point in the Eulerian mesh is not known because the material flows through the mesh from time step to time step. This is the bespoken convective effect of the fluid in the Eulerian framework. For water, since the water's current stress is only dependent on its current strain rate, there is no need to consider the stress history. But for the soil material to be studied in the Eulerian framework, the stress history of the fixed-mesh in the Eulerian framework must be tracked accurately. This is realized in the step of calculating the non-local stress.

To consider the stress history, when the non-local stress is calculated, the field point's coordinate is not based on its current value $x_i^t$. Instead, it is calculated based on the updated coordinate $x_i^{t+dt} = x_i^t + u_i dt$, where $u_i$ is the velocity and $dt$ is the time increment. This treatment means that the current $x_i^t$ field point's material point flows to the new field point $x_i^{t+dt}$ such that the relative distance between this field point and the source point should be calculated based on the update coordinate of the filed point: $d^{mn} = \left| x_i^m - x_i^{n,t+dt} \right|$, where $d^{mn}$ is the relative distance between the source point $x_i^m$ at time $t$ and the field point $x_i^{n,t+dt}$ at time $t+dt$. The relative distance $d^{mn}$ is used to evaluate the weight function $w(d^{mn})$ such that the non-local stress $\hat{\sigma}_{ij}$ is calculated based on the original stress field $\sigma_{ij}$ from all the neighbouring field points:

$$\hat{\sigma}_{ij}^m = \frac{\sum_{n=1} w(d^{mn}) \sigma_{ij}^n}{\sum_{n=1} w(d^{mn})} \tag{1-1}$$

The multiphase material (fluid or solid) structure interaction problem is focused in this paper. To consider the multiphase materials, the interface between different phases of materials must be tracked. The modified level set method is adopted (Touré and Soulaïmani, 2016). This modified level set method adds one more penalty-like diffusive term in the original level set equation such that there is no need to re-initiate the level set variable as required in other level set methods.

Regarding on the interaction algorithm between the fluid and the structure, the idea of the cutting finite element and the immersed finite element is adopted (Zhang *et al*, 2004). The structure's Lagrangian mesh does not need to be compatible with multiphase materials' Eulerian meshes. The cutting points between the structure boundary and the Eulerian mesh edges are calculated. The coupling condition between the structure and the materials is satisfied on these cutting points. The Lagrangian multiplier method is used to reinforce the velocity continuity condition between the structure's boundary and the multiphase materials.

To test the accuracy of this method, various problems are studied, including the fully-buried plate-soil interaction problem, the fully-buried pipe-soil interaction problem and the pipe-water-soil interaction problem with the deformed free surface between the water phase and the soil phase, see Figure 1-1. The elasplastic soil material is regarded as the stress-history-dependent fluid material governed by the N-S equation. The soil material is elastic before the critical yielding point. After yielding point, it is plastic. Moreover, the accumulated plastic strain degenerates the operative soil strength, which is the softening effect. Meanwhile, the strain rate also has some hardening effect. The details will be discussed later.

The literature review for the study of pipe-soil interaction problem is shown here. In (Liu and Marsden, 2018), the multi-scale finite element has been successfully applied to study both solid and fluid materials. The fluid structure interaction problem is also investigated.

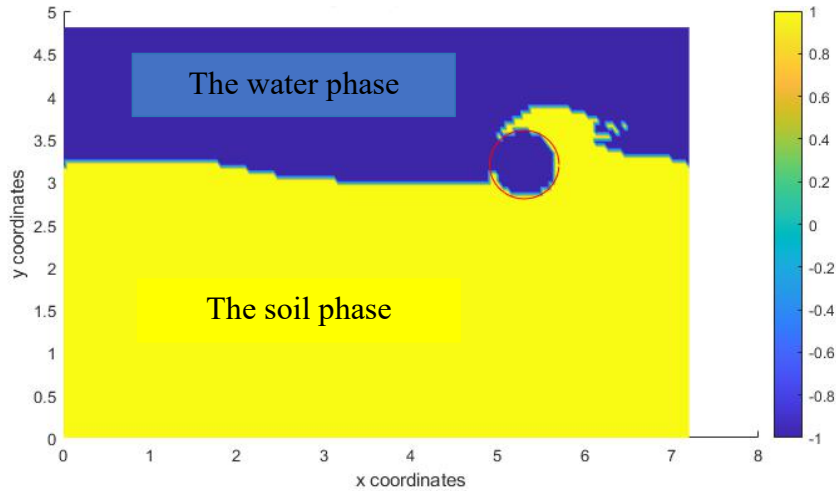

**Figure 1-1** The illustrative figure of the pipe-soil-water interaction analysis

In (Merifield *et al*, 2009), the 2D ALE (Arbitrary Lagrangian and Eulerian) method is applied to study the problem of the pipe's vertical penetration into the soil foundation. The ALE method re-meshes the soil domain based on the movement of pipe's boundary. Moderate deformation problems can be studied well.

A re-mesh based Lagrangian finite element method known as RITSS is proposed by (Wang *et al,* 2009) and (Wang *et al*, 2010). In these papers, the large deformation problem is solved through running a series of small strain linear finite element analyses. The stress mapping between the old and new mesh is realized based on interpolation.

The explicit dynamics software LS-DYNA is applied to simulate the 3D problem of pipe-soil interaction in (Konuk and Yu, 2007) and (Yu and Konuk, 2007), as shown in Figure 1-2. The ALE method is used for the update of mesh. In the explicit dynamics, only the mass matrix is used to calculate the displacement increment without involving the tangential stiffness matrix. As a result, the time step must very small. However, for these transient impact problems, the explicit algorithm is very suitable. This is because the time step must be naturally very small due to the significant dynamical effect. Even in the implicit dynamics algorithm, the time step increment cannot be any larger. In this case, the explicit dynamics becomes more suitable. This phenomena is widely-observed. For example, the pipe dynamical buckling propagation problem analysed by (Li , 2017) is extremely suitable using the explicit dynamics. The time increment is of order -8 (s). The implicit dynamics analysis's time increment is also of the same order in ANSYS LS-DYNA.

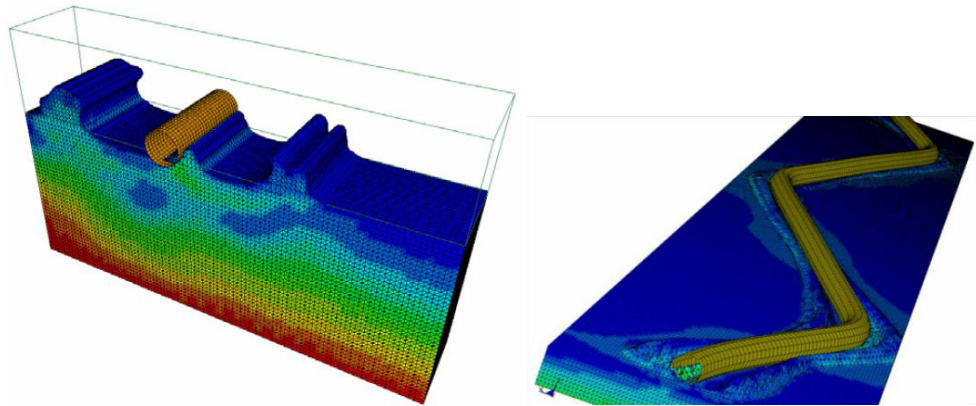

**Figure 1-2** Arbitrary Lagrangian-Eulerian (ALE) model for pipe-soil interaction (Yu and Konuk, 2007) **(a)** 3D plane strain **(b)** Full 3D

In (Dutta *et al*, 2012) and (Dutta *et al* ,2015), an ABAQUS CEL (coupled Eulerian and Lagrangian) analysis is presented for the study of the dynamical embedment processes between the 2D pipe segment and the deformeable soil foundation. The CEL method automatically determines the usage of the Lagrangian framework or Eulerian framework based on the deformation gradient magnitude. The advantage of CEL is that it can handle the large deformation analysis very well. The drawback is that there is only one number to quantify the volume fraction ratio within each cell of the Eulerian mesh. It is not possible to recover the interface exactly unless a very fine mesh is used. The simulation efficiency is not very good.

The DEM (discrete element method) is applied to study the interaction between the 2D pipe segment and the sandy soil (Macaro, 2015), see Figure 1-3. The DEM is a particle-based method. It is Lagrangian. However, since the DEL does not need a mesh, the mesh distortion problem does not arise here. The solution domain is discretized by a cluster of discrete particles. The interacting forces between different particles is modelled using a nonlinear contact laws. The main drawback of the DEM is the stability problem. Meanwhile, an extremely small time step is also required leading to potentially unacceptable computational cost for realistic 3D engineering practice. Moreover, the application of the boundary conditions in these particle-based methods is as convenient as that in those mesh-based methods.

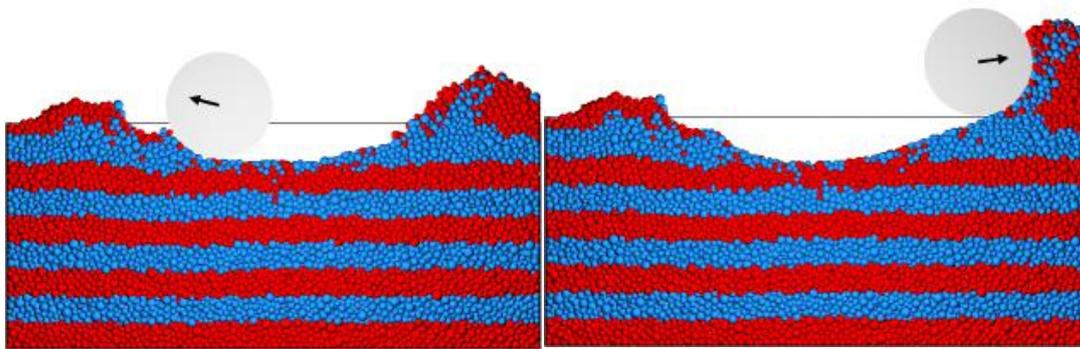

**Figure 1-3** Distinct element method (DEM) model for pipe-soil interaction (Macaro, 2015) (a) Cycle 3, u = 0.5D (b) Cycle 3, u = 2.5D

The SPH (Smooth Particle Hydrodynamics) method is applied to study the complex soil mechanics problem in (Bui and Fukagawa, 2013). The essential components of the SPH theory are presented with enough details, including the definition of the Kernel function, the implementation details of the SPH, the artificial compressibility etc. A velocity-based plastic constitutive formulation is presented.

The outline of the paper is given here. In section 2, the formulation of the variational multi-scale finite element method is derived based on the differential operator. In section 3, the level set method is applied to track the development of the interfaces between different phases of fluids. In section 4, the algorithm of the fluid structure interaction is discussed. In section 5, the history-dependent material constitutive model in the Eulerian (fixed-mesh) framework is discussed. In section 6, the constitutive model for the elasplastic soil material is discussed in the Eulerian framework. The strain-rate-dependent hardening effect and the strain-dependent softening effect are both considered. In section 7, various numerical problems are studied using the proposed method; the numerical results are verified with analytical solutions or numerical results from the ABAQUS CEL method. The problem of the pipe-soil-water interaction is studied using the proposed method.

## 2. VMS-FEM

In this section, the N-S equation with an extra stress term is presented. The corresponding VMS-FEM (Hughes *et al*, 1998) formulation is derived based on the differential operator of the N-S equation. The multiphase fluid case is focused, see Figure 2-1.

The incompressible N-S equation reads:

$$R^M(u_i, p) = \rho \dot{u}_i + \rho u_j u_{i,j} + p_{,i} - \mu u_{i,kk} + \sigma_{ij,j} - b_i = 0$$
$$R^m(u_j) = u_{j,j} = 0 \qquad \text{in} \quad \Omega_I \times t \qquad (2\text{-}1a, b)$$

$$u_i = \bar{u}_i \qquad \text{in} \quad \partial_u \Omega_I \times t \qquad (2\text{-}1c)$$

$$T_{ij} = -p\delta_{ij} + \mu(u_{i,j} + u_{j,i}) + \sigma_{ij} = h_{ij} \qquad \text{in} \quad \partial_h \Omega_I \times t \qquad (2\text{-}1d)$$

$$u_i(t=0) = u_i^0$$
$$p(t=0) = p^0 \qquad \text{in} \quad \Omega_I \times (t=0) \qquad (2\text{-}1e, f)$$

where $u_i$ is the velocity $i$ component, $p$ is the pressure, $\rho$ is density, $\mu$ is the fluid viscosity, $b_i$ is the body force $i$ component, $\sigma_{ij}$ is the extra stress tensor, $\Omega_I$ is the domain of phase $I$, $\partial_u \Omega_I$ is the velocity boundary condition for the domain of phase $I$, $\partial_h \Omega_I$ is the traction boundary condition for the domain of phase $I$, $T_{ij}$ is the total stress tensor, $h_{ij}$ is the boundary traction tensor, $u_i^0$ and $p^0$ are the initial condition for velocity and pressure, $t$ is time.

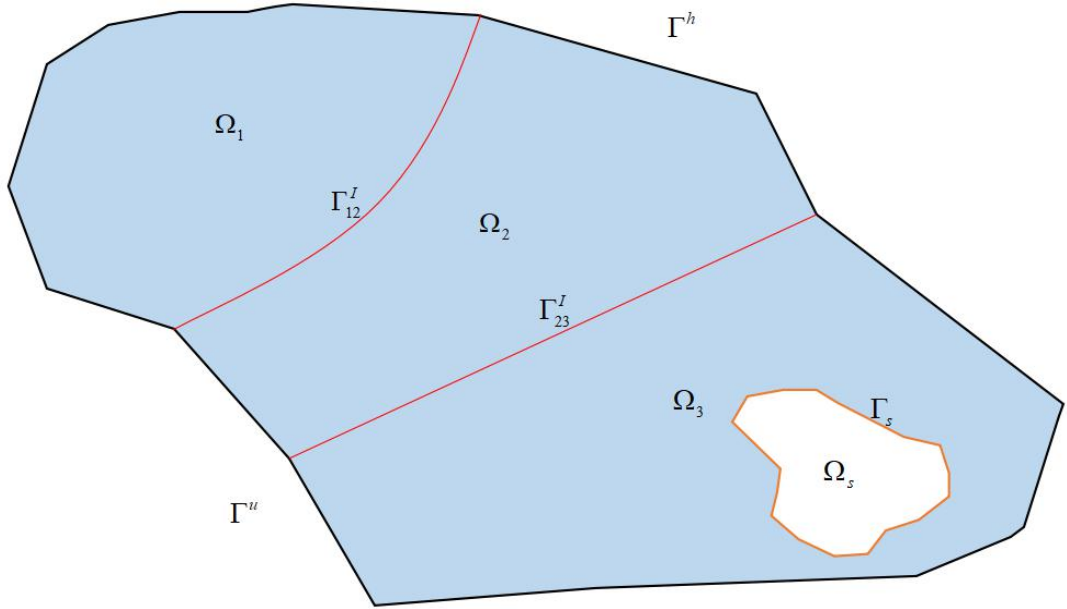

**Figure 2-1** The multiphase fluid structure interaction problem ($\Omega_s$ is the structure domain, $\Gamma_s$ is the structure boundary, $\Omega_i \mid i=1,2,3$ is the fluid domain for fluid material $i$, $\Gamma_{ij}^I$ is the interface between fluid material $i$ and fluid material $j$, $\Gamma^u = \partial_u \Omega$ is the velocity boundary of fluid domain, $\Gamma^h = \partial_h \Omega$ is the velocity boundary of fluid domain)

In this paper, the extra stress tensor $\sigma_{ij}$ represents all the solid-like stress field, which depends on the current deformation gradient and the previous stress history of the current material point. The extra stress tensor $\sigma_{ij}$ is not served as an additional independent degree-of-freedom. It is calculated as a function of velocity $u_i$, which can be regarded as a source term. The calculation of the extra stress $\sigma_{ij}$ will be discussed in Section 5 and 6.

The N-S equation (2-1) can be re-written in the differential operator form:

$$\begin{bmatrix} R^M(u_i, p) \\ R^m(u) \end{bmatrix} = \begin{bmatrix} L^M & L_i^m \\ L_i^m & 0 \end{bmatrix} \begin{bmatrix} u_i \\ p \end{bmatrix} - \begin{bmatrix} f_i \\ 0 \end{bmatrix} = \begin{bmatrix} 0 \\ 0 \end{bmatrix} \qquad (2\text{-}2)$$

where the differential operators are defined as:

$$L^M = L_1^M + L_2^M = \rho\frac{\partial}{\partial t} + \rho u_j \frac{\partial}{\partial x_j} - \mu\frac{\partial^2}{\partial x_k \partial x_k}$$

$$L_j^m = \frac{\partial}{\partial x_j}$$

$$L_1^M = \rho\frac{\partial}{\partial t} + \rho u_j \frac{\partial}{\partial x_j} \qquad (2\text{-}3\text{a, b, c, d, e})$$

$$L_2^M = -\mu\frac{\partial^2}{\partial x_k \partial x_k} = -\left(L_3^M\right)^2 = -\left(\sqrt{\mu}\frac{\partial}{\partial x_k}\right)^2$$

$$f_i = -\sigma_{ij,j} + b_i$$

It is noted that the differential operator form is more convenient to derive the VMS-FEM formulation where the inverse of the operator are used to determine the stability parameters (Hughes *et al*, 2004).

In the VMS-FEM, the main variables are divided into two scales:

$$u_i = \bar{u}_i + \tilde{u}_i$$
$$p = \bar{p} + \tilde{p} \qquad (2\text{-}4\text{a, b})$$

where $\bar{u}_i$ and $\bar{p}$ are the macro velocity and pressure, $\tilde{u}_i$ and $\tilde{p}$ are the micro velocity and pressure.

Thus, (2-2) becomes:

$$\begin{aligned}
\begin{bmatrix} R^M(\bar{u}_i + \tilde{u}_i, \bar{p}+\tilde{p}) \\ R^m(\bar{u}_i + \tilde{u}_i) \end{bmatrix} &= \begin{bmatrix} L^M & L_i^m \\ L_i^m & 0 \end{bmatrix}\begin{bmatrix} \bar{u}_i + \tilde{u}_i \\ \bar{p}+\tilde{p} \end{bmatrix} - \begin{bmatrix} f_i \\ 0 \end{bmatrix} \\
&= \begin{bmatrix} L_1^M + L_2^M & L_i^m \\ L_i^m & 0 \end{bmatrix}\begin{bmatrix} \bar{u}_i \\ \bar{p} \end{bmatrix} + \begin{bmatrix} L_1^M + L_2^M & L_i^m \\ L_i^m & 0 \end{bmatrix}\begin{bmatrix} \tilde{u}_i \\ \tilde{p} \end{bmatrix} - \begin{bmatrix} f_i \\ 0 \end{bmatrix} \\
&= \begin{bmatrix} L_1^M & 0 \\ L_i^m & 0 \end{bmatrix}\begin{bmatrix} \bar{u}_i \\ \bar{p} \end{bmatrix} + \begin{bmatrix} -(L_3^M)^2 & L_i^m \\ 0 & 0 \end{bmatrix}\begin{bmatrix} \bar{u}_i \\ \bar{p} \end{bmatrix} + \begin{bmatrix} L_1^M + L_2^M & L_i^m \\ L_i^m & 0 \end{bmatrix}\begin{bmatrix} \tilde{u}_i \\ \tilde{p} \end{bmatrix} - \begin{bmatrix} f_i \\ 0 \end{bmatrix} = \begin{bmatrix} 0 \\ 0 \end{bmatrix}
\end{aligned} \qquad (2\text{-}5)$$

Set the testing function $v_i, w$ for the momentum equation and the continuity equation, using the similar scale separation rules as in (2-4a, b) for $v_i, w$, we have:

$$v_i = \bar{v}_i + \tilde{v}_i$$
$$w = \bar{w} + \tilde{w} \qquad (2\text{-}6\text{a, b})$$

where $\bar{v}_i, \bar{w}$ are the macro scale testing function for the momentum equation and the continuity equation, and $\tilde{v}_i, \tilde{w}$ are micro scale testing function for the momentum equation and the continuity equation.

Thus, for the N-S equation (2-5), after using (2-6a, b), we have:

$$\sum_I \left\langle \begin{bmatrix} \bar{v}_i + \tilde{v}_i \\ \bar{w} + \tilde{w} \end{bmatrix}^T \begin{bmatrix} R^M(\bar{u}_i + \tilde{u}_i, \bar{p} + \tilde{p}) \\ R^m(\bar{u}_i + \tilde{u}_i) \end{bmatrix} \right\rangle_{\Omega_I}$$

$$= \sum_I \left\langle \begin{bmatrix} \bar{v}_i + \tilde{v}_i \\ \bar{w} + \tilde{w} \end{bmatrix}^T \begin{bmatrix} L^M & L_i^m \\ L_i^m & 0 \end{bmatrix} \begin{bmatrix} \bar{u}_i + \tilde{u}_i \\ \bar{p} + \tilde{p} \end{bmatrix} \right\rangle_{\Omega_I} - \sum_I \left\langle \begin{bmatrix} \bar{v}_i + \tilde{v}_i \\ \bar{w} + \tilde{w} \end{bmatrix}^T \begin{bmatrix} f_i \\ 0 \end{bmatrix} \right\rangle_{\Omega_I}$$

$$= \sum_I \left\langle \begin{bmatrix} \bar{v}_i \\ \bar{w} \end{bmatrix}^T \begin{bmatrix} L^M & L_i^m \\ L_i^m & 0 \end{bmatrix} \begin{bmatrix} \bar{u}_i + \tilde{u}_i \\ \bar{p} + \tilde{p} \end{bmatrix} \right\rangle_{\Omega_I} - \sum_I \left\langle \begin{bmatrix} \bar{v}_i \\ \bar{w} \end{bmatrix}^T \begin{bmatrix} f_i \\ 0 \end{bmatrix} \right\rangle_{\Omega_I} \quad (2\text{-}7)$$

$$+ \sum_I \left\langle \begin{bmatrix} \tilde{v}_i \\ \tilde{w} \end{bmatrix}^T \begin{bmatrix} L^M & L_i^m \\ L_i^m & 0 \end{bmatrix} \begin{bmatrix} \bar{u}_i + \tilde{u}_i \\ \bar{p} + \tilde{p} \end{bmatrix} \right\rangle_{\Omega_I} - \sum_I \left\langle \begin{bmatrix} \tilde{v}_i \\ \tilde{w} \end{bmatrix}^T \begin{bmatrix} f_i \\ 0 \end{bmatrix} \right\rangle_{\Omega_I} = \begin{bmatrix} 0 \\ 0 \end{bmatrix}$$

The equation (2-7) has two scales. The macro scale and the micro scale are assumed to be independent. Thus, for the micro-scale-testing-function ($\tilde{v}_i, \tilde{w}$) related part in (2-7), we have:

$$\sum_I \left\langle \begin{bmatrix} \tilde{v}_i \\ \tilde{w} \end{bmatrix}^T \begin{bmatrix} R^M(\bar{u}_i + \tilde{u}_i, \bar{p} + \tilde{p}) \\ R^m(\bar{u}_i + \tilde{u}_i) \end{bmatrix} \right\rangle_{\Omega_I}$$

$$= \sum_I \left\langle \begin{bmatrix} \tilde{v}_i \\ \tilde{w} \end{bmatrix}^T \left( \begin{bmatrix} L^M & L_i^m \\ L_i^m & 0 \end{bmatrix} \begin{bmatrix} \bar{u}_i \\ \bar{p} \end{bmatrix} - \begin{bmatrix} f_i \\ 0 \end{bmatrix} \right) \right\rangle_{\Omega_I} + \sum_I \left\langle \begin{bmatrix} \tilde{v}_i \\ \tilde{w} \end{bmatrix}^T \begin{bmatrix} L^M & L_i^m \\ L_i^m & 0 \end{bmatrix} \begin{bmatrix} \tilde{u}_i \\ \tilde{p} \end{bmatrix} \right\rangle_{\Omega_I} \quad (2\text{-}8)$$

$$= \sum_I \left\langle \begin{bmatrix} \tilde{v}_i \\ \tilde{w} \end{bmatrix}^T \left( \begin{bmatrix} R^M(\bar{u}_i, \bar{p}) \\ R^m(\bar{u}_i) \end{bmatrix} \right) \right\rangle_{\Omega_I} + \sum_I \left\langle \begin{bmatrix} \tilde{v}_i \\ \tilde{w} \end{bmatrix}^T \begin{bmatrix} L^M & L_i^m \\ L_i^m & 0 \end{bmatrix} \begin{bmatrix} \tilde{u}_i \\ \tilde{p} \end{bmatrix} \right\rangle_{\Omega_I} = \begin{bmatrix} 0 \\ 0 \end{bmatrix}$$

It is expected to determine the micro scale variables $\tilde{u}_i, \tilde{p}$ from (2-8). However, the inverse of equation (2-8) has no explicit formulation, so approximate methods must be used. It is assumed that the micro scale variables can be written as:

$$\begin{bmatrix} \tilde{u}_i \\ \tilde{p} \end{bmatrix} = -\begin{bmatrix} \tau_M & 0 \\ 0 & \tau_m \end{bmatrix} \begin{bmatrix} R^M(\bar{u}_i, \bar{p}) \\ R^m(\bar{u}_i, \bar{p}) \end{bmatrix} \quad (2\text{-}9)$$

where $\tau_M, \tau_m$ are the stability parameter matrices to be obtained from (2-8). The Fourier transformation is adopted to evaluate the values of stability parameters.

For the macro-scale-testing-function related part from (2-7), after using the operator defined in (2-3), we have:

$$\begin{aligned}
&\sum_I \left\langle \begin{bmatrix} \bar{v}_i \\ \bar{w} \end{bmatrix}^T \begin{bmatrix} L^M & L_i^m \\ L_i^m & 0 \end{bmatrix} \begin{bmatrix} \bar{u}_i + \tilde{u}_i \\ \bar{p} + \tilde{p} \end{bmatrix} \right\rangle_{\Omega_I} - \sum_I \left\langle \begin{bmatrix} \bar{v}_i \\ \bar{w} \end{bmatrix}^T \begin{bmatrix} f_i \\ 0 \end{bmatrix} \right\rangle_{\Omega_I} \\
&= \sum_I \left\langle \begin{bmatrix} \bar{v}_i \\ \bar{w} \end{bmatrix}^T \begin{bmatrix} L_1^M + L_2^M & L_i^m \\ L_i^m & 0 \end{bmatrix} \begin{bmatrix} \bar{u}_i \\ \bar{p} \end{bmatrix} \right\rangle_{\Omega_I} + \\
&\sum_I \left\langle \begin{bmatrix} \bar{v}_i \\ \bar{w} \end{bmatrix}^T \begin{bmatrix} L_1^M + L_2^M & L_i^m \\ L_i^m & 0 \end{bmatrix} \begin{bmatrix} \tilde{u}_i \\ \tilde{p} \end{bmatrix} \right\rangle_{\Omega_I} - \sum_I \left\langle \begin{bmatrix} \bar{v}_i \\ \bar{w} \end{bmatrix}^T \begin{bmatrix} f_i \\ 0 \end{bmatrix} \right\rangle_{\Omega_I} \\
&= \sum_I \left\langle \begin{bmatrix} \bar{v}_i \\ \bar{w} \end{bmatrix}^T \begin{bmatrix} L_1^M & 0 \\ L_i^m & 0 \end{bmatrix} \begin{bmatrix} \bar{u}_i \\ \bar{p} \end{bmatrix} \right\rangle_{\Omega_I} + \\
&\sum_I \left\langle \begin{bmatrix} \bar{v}_i \\ \bar{w} \end{bmatrix}^T \begin{bmatrix} -(L_3^M)^2 & L_i^m \\ 0 & 0 \end{bmatrix} \begin{bmatrix} \bar{u}_i \\ \bar{p} \end{bmatrix} \right\rangle_{\Omega_I} - \sum_I \left\langle \begin{bmatrix} \bar{v}_i \\ \bar{w} \end{bmatrix}^T \begin{bmatrix} f_i \\ 0 \end{bmatrix} \right\rangle_{\Omega_I} \\
&+ \sum_I \left\langle \begin{bmatrix} \bar{v}_i \\ \bar{w} \end{bmatrix}^T \begin{bmatrix} L_1^M + L_2^M & L_i^m \\ L_i^m & 0 \end{bmatrix} \begin{bmatrix} \tilde{u}_i \\ \tilde{p} \end{bmatrix} \right\rangle_{\Omega_I} = \begin{bmatrix} 0 \\ 0 \end{bmatrix}
\end{aligned} \quad (2\text{-}10)$$

It is noted that:

$$\begin{bmatrix} -(L_3^M)^2 & L_i^m \\ 0 & 0 \end{bmatrix} = \begin{bmatrix} -L_3^M & L_i^m \\ 0 & 0 \end{bmatrix} \begin{bmatrix} L_3^M & 0 \\ 0 & 1 \end{bmatrix} \quad (2\text{-}11)$$

Thus, using (2-11) and applying integration by part, equation (2-10) becomes:

$$\sum_{I}\left\langle\begin{bmatrix}\bar{v}_i\\\bar{w}\end{bmatrix}^T\begin{bmatrix}L_1^M & 0\\L_i^m & 0\end{bmatrix}\begin{bmatrix}\bar{u}_i\\\bar{p}\end{bmatrix}\right\rangle_{\Omega_I}+\sum_{I}\left\langle\begin{bmatrix}\bar{v}_i\\\bar{w}\end{bmatrix}^T\begin{bmatrix}-L_3^M & L_i^m\\0 & 0\end{bmatrix}\begin{bmatrix}L_3^M & 0\\0 & 1\end{bmatrix}\begin{bmatrix}\bar{u}_i\\\bar{p}\end{bmatrix}\right\rangle_{\Omega_I}$$

$$-\sum_{I}\left\langle\begin{bmatrix}\bar{v}_i\\\bar{w}\end{bmatrix}^T\begin{bmatrix}f_i\\0\end{bmatrix}\right\rangle_{\Omega_I}+\sum_{I}\left\langle\begin{bmatrix}\bar{v}_i\\\bar{w}\end{bmatrix}^T\begin{bmatrix}L_1^M+L_2^M & L_i^m\\L_i^m & 0\end{bmatrix}\begin{bmatrix}\tilde{u}_i\\\tilde{p}\end{bmatrix}\right\rangle_{\Omega_I}$$

$$=\sum_{I}\left\langle\begin{bmatrix}\bar{v}_i\\\bar{w}\end{bmatrix}^T\begin{bmatrix}L_1^M & 0\\L_i^m & 0\end{bmatrix}\begin{bmatrix}\bar{u}_i\\\bar{p}\end{bmatrix}\right\rangle_{\Omega_I}-\sum_{I}\left\langle\begin{bmatrix}-L_3^{M*} & L_i^{m*}\\0 & 0\end{bmatrix}\begin{bmatrix}\bar{v}_i\\\bar{w}\end{bmatrix}^T\begin{bmatrix}L_3^M & 0\\0 & 1\end{bmatrix}\begin{bmatrix}\bar{u}_i\\\bar{p}\end{bmatrix}\right\rangle_{\Omega_I}$$

$$-\sum_{I}\left\langle\begin{bmatrix}\bar{v}_i\\\bar{w}\end{bmatrix}^T\begin{bmatrix}f_i\\0\end{bmatrix}\right\rangle_{\Omega_I}+\sum_{I}\left\langle\begin{bmatrix}\bar{v}_i\\\bar{w}\end{bmatrix}^T\begin{bmatrix}L_1^M+L_2^M & L_i^m\\L_i^m & 0\end{bmatrix}\begin{bmatrix}\tilde{u}_i\\\tilde{p}\end{bmatrix}\right\rangle_{\Omega_I}=\begin{bmatrix}0\\0\end{bmatrix}$$

(2-12)

where the * represents the self-adjoint operator.

The micro scale part variables defined in (2-9) is plugged into (2-12). Then, equation (2-12) becomes:

$$\sum_{I}\left\langle\begin{bmatrix}\bar{v}_i\\\bar{w}\end{bmatrix}^T\begin{bmatrix}L_1^M & 0\\L_i^m & 0\end{bmatrix}\begin{bmatrix}\bar{u}_i\\\bar{p}\end{bmatrix}\right\rangle_{\Omega_I}-\sum_{I}\left\langle\begin{bmatrix}-L_3^{M*} & L_i^{m*}\\0 & 0\end{bmatrix}\begin{bmatrix}\bar{v}_i\\\bar{w}\end{bmatrix}^T\begin{bmatrix}L_3^M & 0\\0 & 1\end{bmatrix}\begin{bmatrix}\bar{u}_i\\\bar{p}\end{bmatrix}\right\rangle_{\Omega_I}-\sum_{I}\left\langle\begin{bmatrix}\bar{v}_i\\\bar{w}\end{bmatrix}^T\begin{bmatrix}f_i\\0\end{bmatrix}\right\rangle_{\Omega_I}$$

$$+\sum_{I}\left\langle\begin{bmatrix}L_1^{M*}+L_2^{M*} & L_i^{m*}\\L_i^{m*} & 0\end{bmatrix}\begin{bmatrix}\bar{v}_i\\\bar{w}\end{bmatrix}^T\begin{bmatrix}\tau_M & 0\\0 & \tau_m\end{bmatrix}\begin{bmatrix}R^M(\bar{u}_i,\bar{p})\\R^m(\bar{u}_i,\bar{p})\end{bmatrix}\right\rangle_{\Omega_I}=\begin{bmatrix}0\\0\end{bmatrix}$$

(2-13)

Using the (2-2), after some re-arrangement, the equation (2-13) becomes:

$$\sum_{I}\left\langle\begin{bmatrix}\bar{v}_i\\\bar{w}\end{bmatrix}^T\begin{bmatrix}L_1^M & 0\\L_i^m & 0\end{bmatrix}\begin{bmatrix}\bar{u}_i\\\bar{p}\end{bmatrix}\right\rangle_{\Omega_I}-\sum_{I}\left\langle\begin{bmatrix}-L_3^{M*} & L_i^{m*}\\0 & 0\end{bmatrix}\begin{bmatrix}\bar{v}_i\\\bar{w}\end{bmatrix}^T\begin{bmatrix}L_3^M & 0\\0 & 1\end{bmatrix}\begin{bmatrix}\bar{u}_i\\\bar{p}\end{bmatrix}\right\rangle_{\Omega_I}$$

$$+\sum_{I}\left\langle\begin{bmatrix}L_1^{M*}+L_2^{M*} & L_i^{m*}\\L_i^{m*} & 0\end{bmatrix}\begin{bmatrix}\bar{v}_i\\\bar{w}\end{bmatrix}^T\begin{bmatrix}\tau_M & 0\\0 & \tau_m\end{bmatrix}\begin{bmatrix}L_1^M+L_2^M & L_i^m\\L_i^m & 0\end{bmatrix}\begin{bmatrix}\bar{u}_i\\\bar{p}\end{bmatrix}\right\rangle_{\Omega_I}$$

(2-14)

$$-\sum_{I}\left\langle\left(\begin{bmatrix}\bar{v}_i\\\bar{w}\end{bmatrix}^T+\begin{bmatrix}L_1^{M*}+L_2^{M*} & L_i^{m*}\\L_i^{m*} & 0\end{bmatrix}\begin{bmatrix}\bar{v}_i\\\bar{w}\end{bmatrix}^T\begin{bmatrix}\tau_M & 0\\0 & \tau_m\end{bmatrix}\right)\begin{bmatrix}f_i\\0\end{bmatrix}\right\rangle_{\Omega_I}=\begin{bmatrix}0\\0\end{bmatrix}$$

In (2-14), the first line is the traditional Galerkin finite element part of the N-S equation, the first term is the inertia and convective term, the second term is the diffusive and pressure gradient term where the integration by part is applied. For the second line, it is the VMS micro-scale-related part. The third line shows all of the force terms.

For the third line of (2-14), using (2-3e) and applying the integration by parts, we have:

$$\sum_I \left\langle \left( \begin{bmatrix} \bar{v}_i \\ \bar{w} \end{bmatrix}^T + \begin{bmatrix} L_1^{M*} + L_2^{M*} & L_i^{m*} \\ L_i^{m*} & 0 \end{bmatrix} \begin{bmatrix} \bar{v}_i \\ \bar{w} \end{bmatrix} \right)^T \begin{bmatrix} \tau_M & 0 \\ 0 & \tau_m \end{bmatrix} \begin{bmatrix} f_i \\ 0 \end{bmatrix} \right\rangle_{\Omega_I}$$

$$= \sum_I \left\langle \frac{\partial}{\partial x_j} \left( \begin{bmatrix} \bar{v}_i \\ \bar{w} \end{bmatrix}^T + \begin{bmatrix} L_1^{M*} + L_2^{M*} & L_i^{m*} \\ L_i^{m*} & 0 \end{bmatrix} \begin{bmatrix} \bar{v}_i \\ \bar{w} \end{bmatrix} \right)^T \begin{bmatrix} \tau_M & 0 \\ 0 & \tau_m \end{bmatrix} \begin{bmatrix} \sigma_{ij} \\ 0 \end{bmatrix} \right\rangle_{\Omega_I}$$

$$+ \sum_I \left\langle \left( \begin{bmatrix} \bar{v}_i \\ \bar{w} \end{bmatrix}^T + \begin{bmatrix} L_1^{M*} + L_2^{M*} & L_i^{m*} \\ L_i^{m*} & 0 \end{bmatrix} \begin{bmatrix} \bar{v}_i \\ \bar{w} \end{bmatrix} \right)^T \begin{bmatrix} \tau_M & 0 \\ 0 & \tau_m \end{bmatrix} \begin{bmatrix} b_i \\ 0 \end{bmatrix} \right\rangle_{\Omega_I}$$

$$\approx \sum_I \left\langle \frac{\partial}{\partial x_j} \begin{bmatrix} \bar{v}_i \\ \bar{w} \end{bmatrix}^T \begin{bmatrix} \sigma_{ij} \\ 0 \end{bmatrix} \right\rangle_{\Omega_I} + \sum_I \left\langle \left( \begin{bmatrix} \bar{v}_i \\ \bar{w} \end{bmatrix}^T + \begin{bmatrix} L_1^{M*} + L_2^{M*} & L_i^{m*} \\ L_i^{m*} & 0 \end{bmatrix} \begin{bmatrix} \bar{v}_i \\ \bar{w} \end{bmatrix} \right)^T \begin{bmatrix} \tau_M & 0 \\ 0 & \tau_m \end{bmatrix} \begin{bmatrix} b_i \\ 0 \end{bmatrix} \right\rangle_{\Omega_I}$$

(2-15)

In this step, the VMS micro-scale-related part is ignored for the extra stress tensor term $\sigma_{ij}$. It is assumed that the solid-like stress only affects the macro scale.

For the update of the extra stress force term in (2-15), it will be discussed in section 5 and 6. This term is handled explicitly.

In order to solve (2-15), the incremental form is used:

$$\begin{bmatrix} K_{uu}^{t+dt} & K_{up}^{t+dt} \\ K_{pu}^{t+dt} & K_{pp}^{t+dt} \end{bmatrix} \begin{bmatrix} d\hat{u} \\ d\hat{p} \end{bmatrix} = \begin{bmatrix} R_t^M \\ R_t^m \end{bmatrix} \quad (2\text{-}16)$$

where $\begin{bmatrix} R_t^M \\ R_t^m \end{bmatrix}$ is the residue of the N-S equation (2-15) from the last time step,

$\begin{bmatrix} K_{uu}^{t+dt} & K_{up}^{t+dt} \\ K_{pu}^{t+dt} & K_{pp}^{t+dt} \end{bmatrix}$ is the tangential stiffness matrix by linearizing equation (2-15) with respect to the main variables $u, p$.

The explicit time integration is adopted:

$$\begin{bmatrix} d\hat{u} \\ d\hat{p} \end{bmatrix} = \begin{bmatrix} K_{uu}^{t+dt}(\hat{u}^t, \hat{p}^t) & K_{up}^{t+dt}(\hat{u}^t, \hat{p}^t) \\ K_{pu}^{t+dt}(\hat{u}^t, \hat{p}^t) & K_{pp}^{t+dt}(\hat{u}^t, \hat{p}^t) \end{bmatrix} \backslash \begin{bmatrix} R_t^M(\hat{u}^t, \hat{p}^t) \\ R_t^m(\hat{u}^t, \hat{p}^t) \end{bmatrix} \quad (2\text{-}17\text{a})$$

where

$$\hat{u}^{t+dt} = \hat{u}^t + d\hat{u}$$
$$\hat{p}^{t+dt} = \hat{p}^t + d\hat{p}$$
(2-17b, c)

Equation (2-17a) is the main equation to evolve the numerical simulation.

## 3. Level set method

In this section, the level set method (Touré and Soulaïmani, 2016) for tracking the interface development between different phases of fluids is presented. The level set method is a widely used interface tracking method in the Eulerian fixed mesh framework.

The original level set equation is:

$$R_\phi = \dot{\phi} + u_i \phi_{,i} = 0 \quad (3\text{-}1)$$

where $\phi$ is the level set, which represents the shortest distance between the current material point to the current interface line between different phases of materials.

For the purpose of simplification, we only consider two phase materials in this section. However, the extension to the multiphase is straightforward by setting $\phi$ as a multi-dimensional vector. The dimension of vector $\phi$ is the number of total types of potential interfaces between different phases of materials. For example, if we have three phases of materials, there can be three types of interfaces.

For the two phase fluid case, the sign of $\phi$ is defined as follow.

For the material 1 $(x_1, x_2) \in \Omega_1$, the sign of the level set $\phi$ is positive:

$$\text{sign}(\phi(x_i)) = 1 \quad (3\text{-}2\text{a})$$

For the material 2 $(x_1, x_2) \in \Omega_2$, the sign of the level set $\phi$ is negative:

$$\text{sign}(\phi(x_i)) = -1 \tag{3-2b}$$

For the interface between material 1 and material 2 $(x_1, x_2) \in \partial\Omega_{12} = \Omega_1 \cap \Omega_2$:

$$\phi(x_i) = 0 \tag{3-2c}$$

where $\partial\Omega_{12}$ is the interface between the material 1 domain $\Omega_1$ and the material 2 domain $\Omega_2$.

The norm of the spatial gradient of the level set must be unite. The following condition for $\phi$ must be satisfied:

$$|\phi_{,i}\phi_{,i}| = 1 \tag{3-3}$$

The norm of the gradient of level set should be naturally always one in mathematics because the level set $\phi$ is the distance. However, if no special treatment is used, with the time involving in equation (3-1), the condition (3-3) may not be satisfied numerically. Some researchers choose to re-initiate the level set $\phi$ after some time steps. In this paper, we adopt another method that can avoid the re-initiation. As reported in (Touré and Soulaïmani, 2016), the following modified level set equation is adopted:

$$\phi_{,t} + u_i \phi_{,i} - (\lambda_1 k_1 \phi_{,j})_{,j} = 0 \tag{3-4}$$

where $k_1 = 1 - \dfrac{1}{|\phi_{,i}\phi_{,i}|}$ and $\lambda_1 = \beta_1 \dfrac{h_{ele}^2 \sqrt{u_i u_i}}{2}$, which adds a penalty-like term to reinforce the constraint condition $|\phi_{,i}\phi_{,i}| = 1$.

For the level set governing equation (3-4), the weak form is:

$$\begin{aligned}&<w, \dot{\phi}> + <w, u_i \phi_{,i}> - <w, (\lambda_1 k_1 \phi_{,j})_{,j}> = \\ &<w, \dot{\phi}> + <w, u_i \phi_{,i}> + <w_{,j}, \lambda_1 k_1 \phi_{,j}> = 0\end{aligned} \tag{3-5}$$

The idea of the VMS-FEM is also applied to the level set equation:

$$\phi = \bar{\phi} + \tilde{\phi} \tag{3-6}$$

where the micro scale level set is defined as proportional to the residue of (3-4):

$$\tilde{\phi} = -\tau_\phi R_\phi = -\tau_\phi \left( \phi_{,t} + u_i \phi_{,i} - (\lambda_1 k_1 \phi_{,j})_{,j} \right) \tag{3-7}$$

For the stability parameter $\tau_\phi$, it is defined in (Touré and Soulaïmani, 2016).

Plunge (3-7) in (3-5), we have:

$$\begin{aligned}
&<w, \dot{\bar{\phi}}> + <w, \dot{\tilde{\phi}}> + <w, u_i \bar{\phi}_{,i}> + <w, u_i \tilde{\phi}_{,i}> \\
&+ <w_{,j}, \lambda_1 k_1 \bar{\phi}_{,j}> + <w_{,j}, \lambda_1 k_1 \tilde{\phi}_{,j}> = 0
\end{aligned} \tag{3-8}$$

Then, the integration by parts is applied to the last term of (3-8). Using (3-7), the weak form of the level set equation (3-8) becomes:

$$\begin{aligned}
&<w, \bar{\phi}_{,t}> + <w, \tilde{\phi}_{,t}> + <w, u_i \bar{\phi}_{,i}> - <u_i w_{,i}, \tilde{\phi}> \\
&+ <w_{,j}, \lambda_1 k_1 \bar{\phi}_{,j}> + <w_{,j}, \lambda_1 k_1 \tilde{\phi}_{,j}> \\
&= <w, \bar{\phi}_{,t}> - <w, \tau_\phi \bar{\phi}_{,tt}> - <w, \tau_\phi u_i \bar{\phi}_{,ti}> \\
&- <w_{,j}, \tau_\phi (\lambda_1 k_1 \bar{\phi}_{,tj})> + <w, u_i \bar{\phi}_{,i}> \\
&+ <u_i w_{,i}, \tau_\phi \bar{\phi}_{,t}> + <u_i w_{,i}, \tau_\phi u_j \bar{\phi}_{,j}> \\
&+ <(u_i w_{,i})_{,j}, \tau_\phi (\lambda_1 k_1 \bar{\phi}_{,j})> + <w_{,j}, \lambda_1 k_1 \bar{\phi}_{,j}> \\
&- <w_{,j}, \lambda_1 k_1 \tau_\phi \bar{\phi}_{,tj}> - <w_{,j}, \lambda_1 k_1 \tau_\phi u_i \bar{\phi}_{,ij}> \\
&- <w_{,jj}, \tau_\phi (\lambda_1 k_1)^2 \bar{\phi}_{,ll}> = 0
\end{aligned} \tag{3-9a}$$

All of the second order derivative terms in (3-9a) are ignored as the linear interpolation function is used for the level set. Then, the weak form (3-9a) becomes:

$$\begin{aligned}
&<w, \bar{\phi}_{,t}> + <w, u_i \bar{\phi}_{,i}> + <u_i w_{,i}, \tau_\phi \bar{\phi}_{,t}> \\
&+ <u_i w_{,i}, \tau_\phi u_j \bar{\phi}_{,j}> + <w_{,j}, \lambda_1 k_1 \bar{\phi}_{,j}> = 0
\end{aligned} \tag{3-9b}$$

The interpolation of the macro scale $\phi$ is:

$$\bar{\phi} = N_\phi \hat{\phi} \tag{3-10}$$

where $N_\phi$ is the basis function matrix for the level set, $\hat{\phi}$ is the degree-of-freedom vector.

Thus, the following matrices are obtained form (3-9b) after using (3-10):

$$M_e^\phi = \int_{V_e} \left(N_\phi + \tau_\phi (u_i N_\phi)_{,i}\right)^T N_\phi dV$$

$$C_e^\phi = \int_{V_e} \left(N_\phi + \tau_\phi (u_i N_\phi)_{,i}\right)^T (u_j N_\phi)_{,j} dV \qquad \text{(3-11a, b, c)}$$

$$K_e^\phi = \int_{V_e} N_{\phi,j}{}^T \lambda_1 k_1 N_{\phi,j} dV$$

Then, the discrete form of (3-9b) is:

$$M_e^\phi \dot{\hat{\phi}} + (C_e^\phi + K_e^\phi)\hat{\phi} = 0 \qquad \text{(3-12)}$$

The explicit time integration is also used. Thus, we have:

$$M_e^\phi \frac{d\hat{\phi}}{dt} + (C_e^\phi + K_e^\phi)(\hat{\phi}^t + d\hat{\phi}) = 0$$
$$\Updownarrow \qquad \text{(3-13a, b)}$$
$$\left(\frac{M_e^\phi}{dt} + (C_e^\phi + K_e^\phi)\right) d\hat{\phi} = -(C_e^\phi + K_e^\phi)\hat{\phi}^t$$

## 4. Fluid structure interaction

The idea of the cutting finite element and the immersed finite element or immersed boundary method (Zhang *et al*, 2004) is adopted to handle the coupling condition between the structure (Lagrangian framework) and the fluid material (Eulerian framework) phases on the fluid structure interface $\Gamma_s$, see Figure 2-1, Figure 4-1 and Figure 4-2. The advantage of this coupling method is that a non-compatible mesh between the structure and the fluid materials can be used such that the regular quadrilateral mesh can be used for the fluid domain in the Eulerian framework, thus allowing a more physics-reasonable definition of the micro scale characteristic length, which is very important in the determination of the stability parameters for the VMS method, see the references in (Hughes *et al*, 2004) for more detail. The solution domain of fluid materials can be defined as a large rectangular domain. The boundary conditions of the fluid domain can be satisfied in the same way of the cutting finite element, see Figure 4-2.

To reinforce the coupling condition between the structure and the fluid, the first step is to find all of the interaction/cutting points between the fluid Eulerian mesh and the boundary of the structure, as shown in Figure 4-1. This is the case of the pipe structure interacting with the fluid domain diecretized by rectangular mesh. As Figure 4-1 shows, the pipe structure (red line) cuts the Eulerian mesh on these cutting points (red circles).

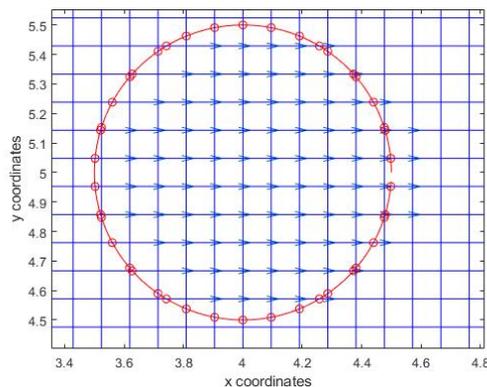

**Figure 4-1** Interaction points between Lagrangian structure (a pipe) and the Eulerian mesh (fluid)

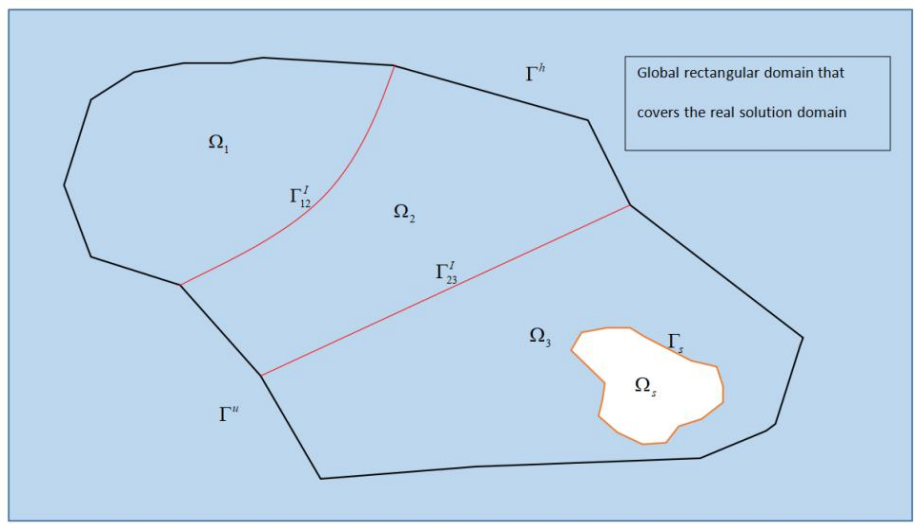

**Figure 4-2** The global rectangular domain that covers the real fluid solution domain

The coupling condition between the structure and the fluid is applied in this way. The nodes inside the current structure domain will be assigned a velocity that is the same as the current velocity of the structure. For the nodes totally outside the current structure domain, the velocity are free and to be solved.

The nodes on the edges that are being cut by the current structure boundary are not totally free. It is important that the cutting point's velocity must be the same with the current velocity of structure. The cutting point velocity is a linear combination of the velocity values of the two nodes of the edge being cut. The coefficients of the linear combination is dependent on the relative distance between the cutting point and the two nodes. Three cases are considered.

**Case 1: non-slip friction model**

In this case, it is assumed that the structure's velocity is totally same with the fluid velocity on the cutting points. The cutting point's velocity is:

$$\vec{u}_{Int} = \alpha \vec{u}_L + (1-\alpha)\vec{u}_R \tag{4-1}$$

where the indices $L$ and $R$ refer to the left and right nodes of the element edge that is being cut by the structure boundary. The coordinate of the cutting point is $P_{Int} = \alpha P_L + (1-\alpha)P_R$, where $P_L, P_R$ are the left and right node coordinates. Obvious, $\alpha \in [0,1]$.

For the non-slip friction model, we have:

$$\vec{u}_{Int} = \alpha \vec{u}_L + (1-\alpha)\vec{u}_R = \vec{u}_s \tag{4-2}$$

where $\vec{u}_s$ is the current velocity of the structure.

After collecting all the constraint conditions (4-2) for all cutting points, we can get the following augmented part for our system:

$$A\hat{u} = b \tag{4-3}$$

Thus, the global system of equations becomes:

$$\begin{bmatrix} K_{uu}^{t+dt} & K_{up}^{t+dt} & A^T \\ K_{pu}^{t+dt} & K_{pp}^{t+dt} & 0 \\ A & 0 & 0 \end{bmatrix} \begin{bmatrix} d\hat{u} \\ d\hat{p} \\ \lambda_s \end{bmatrix} = \begin{bmatrix} R_m^t \\ R_c^t \\ b - A\hat{u}^t \end{bmatrix} \tag{4-4}$$

where $\lambda_s$ is the Lagrangian multiplier corresponding to the constraint condition $A\hat{u} = b$. $\lambda_s$ is the reaction force between the structure and the fluid.

**Case 2: slip friction model**

In this case, only the velocity along the normal direction of the structure boundary $\partial \Omega_s$ is coupled between the structure and the Eulerian fluid:

$$\vec{u}_{Int} \cdot \vec{n}_{Int} = \left( \alpha \vec{u}_L + (1-\alpha)\vec{u}_R \right) \cdot \vec{n}_{Int} = \vec{u}_s \cdot \vec{n}_{Int} \tag{4-5}$$

where $\vec{n}_{Int}$ is the normal vector of the interaction point on $\partial \Omega_s$.

These constraint conditions are also written in matrix form after collecting for all the cutting points:

$$A_s \hat{u} = b_s \tag{4-6}$$

Thus, the global system becomes:

$$\begin{bmatrix} K_{uu}^{t+dt} & K_{up}^{t+dt} & A_s^T \\ K_{pu}^{t+dt} & K_{pp}^{t+dt} & 0 \\ A_s & 0 & 0 \end{bmatrix} \begin{bmatrix} d\hat{u} \\ d\hat{p} \\ \lambda_s \end{bmatrix} = \begin{bmatrix} R_m^t \\ R_c^t \\ b_s - A_s \hat{u}^t \end{bmatrix} \tag{4-7}$$

Note that in this case, the Lagrangian multiplier $\lambda_s$ only includes the reaction force in the normal direction of $\partial \Omega_s$.

**Case 3: general frictional model**

The general frictional model is an intermediate model between the previous two models. The coupling condition between the structure and the fluid is not as strong as Model 1: non-slip model. But the coupling condition is stronger than Model 2: slip model. In this case, we still only constrain the normal velocity between the structure and fluid, which is same as equation (4-5) in Model 2. The difference is that, once we have the solution of $\lambda_s$ (the normal force between solid and fluid) from (4-7), the tangential friction force are calculated and added into the system as one more force:

$$f_{fric} = -sign(u_t)\mu_{firc}\lambda_s, u_t = (\vec{u}_{Int} - \vec{u}_s)\cdot \vec{t} \tag{4-8}$$

where $\vec{t}$ is the tangential unit vector of the cutting point on $\partial\Omega_s$, $u_t$ is the relative tangential velocity between solid and fluid, and $\mu_{fric}$ is the friction coefficient.

Thus, the global system for Model 3 is:

$$\begin{bmatrix} K_{uu}^{t+dt} & K_{up}^{t+dt} & A_s^T \\ K_{pu}^{t+dt} & K_{pp}^{t+dt} & 0 \\ A_s & 0 & 0 \end{bmatrix} \begin{bmatrix} d\hat{u} \\ d\hat{p} \\ \lambda_s \end{bmatrix} = \begin{bmatrix} f_{firc} + R_m^t \\ R_c^t \\ b_s - A_s \hat{u}^t \end{bmatrix} \tag{4-9}$$

In (4-9), one more force term coming from (4-8) is added into the residue of the first line on the right hand side.

## 5. Stress-history-dependent materials in Eulerian framework

In this section, a novel method to study the history-dependent constitutive behaviour in the Eulerian framework is presented, see Figure 5-1. History-dependent constitutive models are widely-used for many solid materials in Lagrangian framework. For the soil material, the elasplastic material model is suitable. The elasplastic constitutive model subroutine needs to know the previous stress history, the previous plastic strain and the previous internal variables such that it can give the increment of stress, plastic strain and internal variable (Simo and Hughes, 2006). This is trivial in Lagrangian framework, where the mesh is always bonded with the materials. However, in the Eulerian framework (fixed mesh), tracking the stress history and the plastic strain history of the material point is not trivial because the material point flows through the fixed mesh from time to time.

In this paper, we propose a novel method to track the stress history of the material point in the Eulerian framework. This is realized in the step of calculating the non-local stress, the non-local plastic strain and the non-local internal variables.

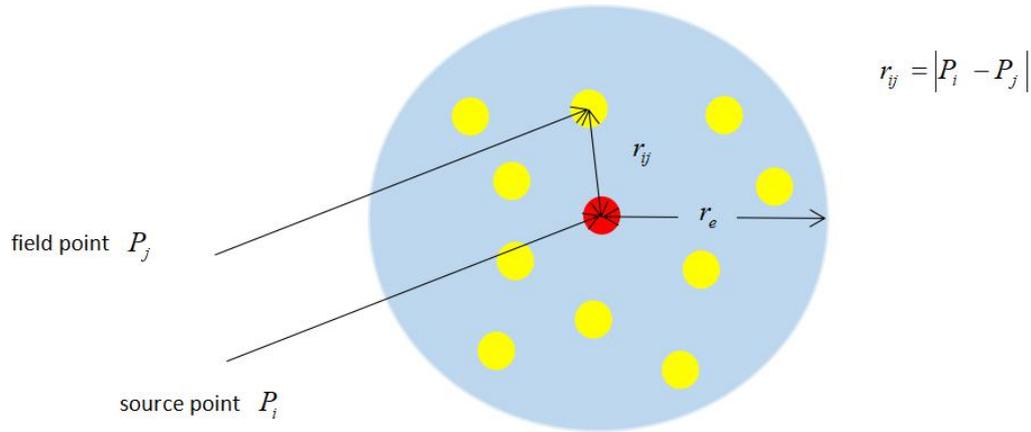

**Figure 5-1 (a)** The traditional non-local stress/strain theory (red circle: source point; yellow circle: field point; blue region: effect radius zone $r < r_e$)

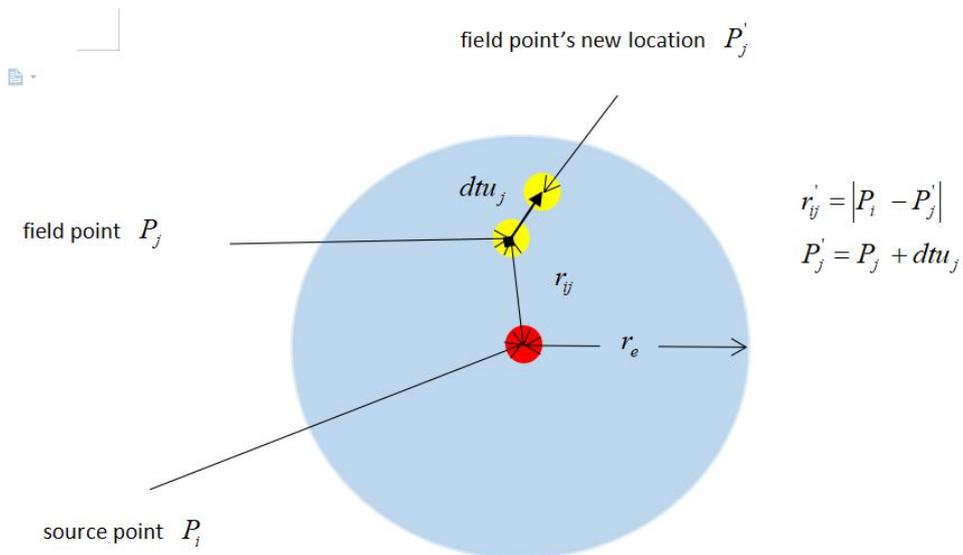

**Figure 5-1 (b)** The proposed non-local stress/strain theory (red circle: source point; yellow circle: field point; blue region: effect radius zone $r < r_e$)

The use of the non-local stress to replace the original stress is very helpful to increase the stability of a numerical method. Highly non-linear problems are always associated with some noise in the numerical solution, especially for the stress and strain fields. To counteract the noise, the non-local stress is used to replace the original stress.

The key idea of the non-local stress/strain theory is to use the values of the stress/strain from all the neighbouring field points to evaluate the stress/strain value of the given target source point. It can be also regarded as an relative-distance-based averaging algorithm to smooth the stress/strain solutions using the Delta function as the weight function (a function about the relative distance between source point and field point).

The non-local stress $\hat{\sigma}_{mn}$ and the non-local strain $\hat{\varepsilon}_{mn}$ are defined as:

$$\hat{\sigma}_{mn} = \frac{\int_{V_r} G(r)\sigma_{mn} dr}{\int_{V_r} G(r) dr}$$

$$\hat{\varepsilon}_{mn} = \frac{\int_{V_r} G(r)\varepsilon_{mn} dr}{\int_{V_r} G(r) dr}$$

(5-1a, b)

where the Delta function is defined as $G(r) = \frac{1}{\sqrt{2\pi}} e^{-\left(\frac{r}{r_0}\right)^2}$, $r$ is the relative distance between the field point and source point, $r_0$ controls the compactness of the Delta function, $V_r$ is the effect domain.

In discrete form, the non-local stress and the non-local strain are defined as:

$$\hat{\sigma}_{mn}^i = \frac{\sum_{j=1} G(r_{ij})\sigma_{mn}^j}{\sum_{j=1} G(r_{ij})}$$

$$\hat{\varepsilon}_{mn}^i = \frac{\sum_{j=1} G(r_{ij})\varepsilon_{mn}^j}{\sum_{j=1} G(r_{ij})}$$

(5-2a, b)

where $r_{ij} = P_i - P_j$ is the relative distance between the source point $P_i$ and the field point $P_j$.

For the current VMS-FEM using the Eulerian framework, we propose the following modified non-local stress/strain formulations:

$$\hat{\sigma}_{mn}^i = \frac{\sum_{j=1} G(r_{ij}^*)\sigma_{mn}^j}{\sum_{j=1} G(r_{ij}^*)}$$

$$\hat{\varepsilon}_{mn}^i = \frac{\sum_{n=1} G(r_{ij}^*)\varepsilon_{mn}^j}{\sum_{j=1} G(r_{ij}^*)}$$

(5-3a, b)

where $r_{ij}^*$ is the updated relative distance between the source point $r_i^*$ and the field point $r_j^*$:

$$r_{ij}^* = r_i^* - r_j^* \tag{5-4a}$$

where

$$\begin{aligned} r_i^* &= r_i + u_i dt \\ r_j^* &= r_j + u_j dt \end{aligned} \tag{5-4b, c}$$

In this case, the convective effect is automatically considered such that we can track the stress/strain history in the Eulerian framework.

## 6. Constitutive model

In the section 5, the non-local stress/strain method is modified to realize the convective effect of the previously stress/strain for the Eulerian framework. In this section, the update algorithm of the history-dependent stress and plastic strain is discussed.

Given the current velocity field from the VMS-FEM, the deformation gradient is:

$$F_{ij} = \frac{\partial u_i}{\partial x_j} dt = \frac{\partial u_i}{\partial \xi_k} \frac{\partial \xi_k}{\partial x_j} dt \tag{6-1}$$

where $\xi_k$ is the natural coordinate, $J_{jk}^{-1} = \frac{\partial \xi_k}{\partial x_j}$ is the inverse of the Jacobian matrix.

Given (6-1), the phenomena plasticity model developed in (Simo and Hughes, 2006) is used to update the stress and plastic strain.

For the soil material, the yielding function is discussed here. The smoothed Tresca criterion is adopted:

$$\psi = \sqrt{0.25(s_{11} - s_{22})^2 + s_{12}^2} - \sigma_y \tag{6-2}$$

where $s_{ij}$ is the deviatoric part of $\sigma_{ij}$, $\sigma_y$ is the Tresca yielding stress.

For the yielding stress, the plastic-strain-dependent softening effect and the strain-rate-dependent hardening effect are both considered. Thus, we define:

$$\sigma_y = \sigma_0 \times \left(\delta_{rem} + (1 - \delta_{rem})e^{-\xi/\xi_{95}}\right) \times \left(1 + \mu_{rate} \log\left(\max(\dot{\gamma}_{max}, \dot{\gamma}_{ref})/\dot{\gamma}_{ref}\right)\right) \tag{6-3}$$

where $\sigma_0$ is the original strength, $\delta_{rem} = S_t^{-1} = 1/3.2$ is the inverse of the soil sensitivity, $\xi$ is the accumulated effect plastic strain, $\xi_{95} = 10$ is the accumulated effect plastic strain under which 95% of soil is remoulding, $\mu_{rate} = 0, 0.05, 0.1$ is the viscosity in the rate-dependent plastic model, $\dot{\gamma}_{max} = \dot{\varepsilon}_1 - \dot{\varepsilon}_2$ is the maximum shear strain rate (principle shear strain rate), $\dot{\varepsilon}_1, \dot{\varepsilon}_2$ are the principle strain rate, $\dot{\gamma}_{ref}$ is reference shear strain rate.

It is noted that the bisection method is used to determine $d\lambda$, the internal variable, in the calculation of the consistency condition, see (Simo and Hughes, 2006). The yielding space defined in (6-3) is also dependent on $d\lambda$, causing nonlinearity in the consistency equation, a scalar nonlinear equation governing the $d\lambda$.

$$\xi = d\lambda + \xi_t, \dot{\gamma}_{max} = d\lambda/dt \tag{6-4a, b}$$

where $\xi_t$ is the previously accumulated plastic effect strain.

## 7. Numerical results

Three problems are studied. In section 7.1, the fully buried plate dragging in the soil domain is analysed. The non-dimensional resistance force is studied and compared with the analytical solution. In section 7.2, the fully buried pipe dragging in the soil domain is studied. In these two problems, the strain-rate-dependent hardening effect and the plastic-strain-dependent softening effect are not considered. Only in this case can we make a comparison with the analytical solutions. The yielding stress is a constant. There is only one phase of soil material.

In section 7.3, the pipe-soil interaction is studied. The vertical penetration and the horizontal movement are both considered. The gravity, the softening effect and the strain-rate-dependent hardening effects are all considered. The numerical results are compared with the SLA method in (Kong, 2015) and the ABAQUS CEL.

In these numerical simulation, the traditional Lagrangian linear interpolation is used for the velocity, pressure and level set fields if not otherwise stated. The bubble function is used to enrich the velocity field in some of the following simulations. The quadratic mesh is used for the fluid domain discretization.

### 7.1 Fully buried plate

In this section, the resistance force of dragging the thin plate fully embedded in the soil domain is studied. The soil constitutive model is the bi-linear elasplastic model based on Tresca shear stress criterion such that the we can compare the numerical solution with the analytical solution reported by (Rowe and Davis, 1982).

According to the analytical solution, the non-dimensional resistance force between the soil and the thin plate is:

$$N_{plate} = \frac{P}{Bs_u} = 11.42 \qquad (7\text{-}1)$$

where $B$ is the width of the plate, $s_u$ is the Tresca undrained shear strength, and $P$ is the resistance force. This analytical solution is based on the rigid-plastic limit analysis and is only valid when the speed of dragging the plate is small enough such that no inertia effect is considered.

In this study, the dragging velocity is $v_d = 0.1B$ m/s. The time increment is 0.001s. The solution domain is $l_a \times l_b = 5B \times 6B$. Considering the symmetrical characteristic of the problem, only half solution domain is used. At $y = 3B$, the symmetrical boundary condition is applied. For the other three edges, the velocity field is fixed.

In Figure 7-1 (a-b), the non-dimensional resistance force is shown. The curve of $N_{plate}$ verses the time step is wave-like. The average value is 12, which is roughly consistent with the analytical solution (7-1). In Figure 7-1 (a), a linear interpolation function is used for all the variables (velocity, pressure and level set). Figure 7-1 (b) shows the case when the bubble function is used to enrich the velocity solution interpolation space. The curve of $N_{plate}$ verses time step is also a wave-like, but the average value is much closer to the analytical solution (11.42). Thus, the bubble function makes the calculation of the resistance force more accurate. In Figure 7-2, the velocity field is shown.

For the simulation in Figure 7-3 and 7-4, the solution domain is larger: $l_a \times l_b = 20B \times 20B$. The average value of the non-dimensional plate-soil resistance force $P/Bs_u$ is smaller than that in Figure 7-1. But the difference is very small. Thus, it is concluded that the solution domain $l_a \times l_b = 5B \times 6B$ in Figure 7-1 (a, b) is acceptable. The boundary effect is not dominating.

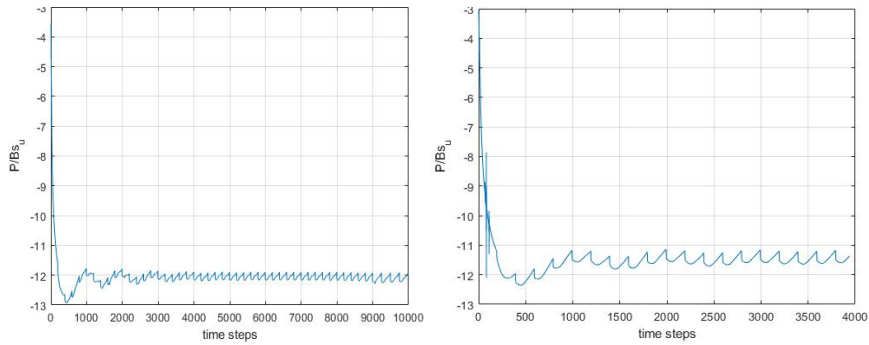

**Figure 7-1 (a, b)** Non-dimensional force without bubble function (left) and with bubble function (right) verses time for thin plate dragged in soil (soil solution domain: $l_a \times l_b = 5B \times 6B$)

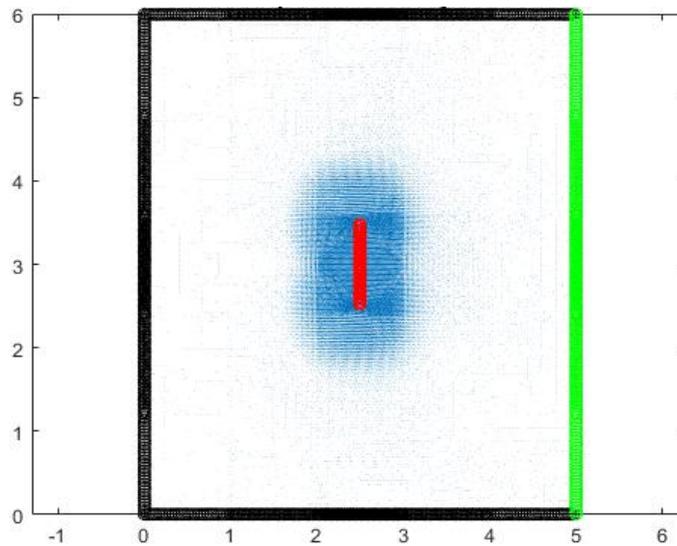

**Figure 7-2** The velocity quiver of soil (with bubble enrichment function) for the problem of thin plate dragged in soil

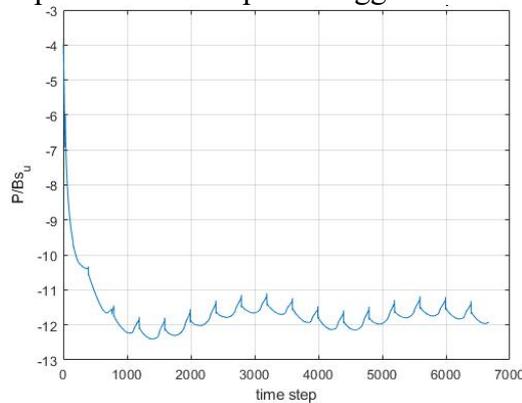

**Figure 7-3** Non-dimensional force verses time for thin plate dragged in soil (soil solution domain: $l_a \times l_b = 20B \times 20B$)

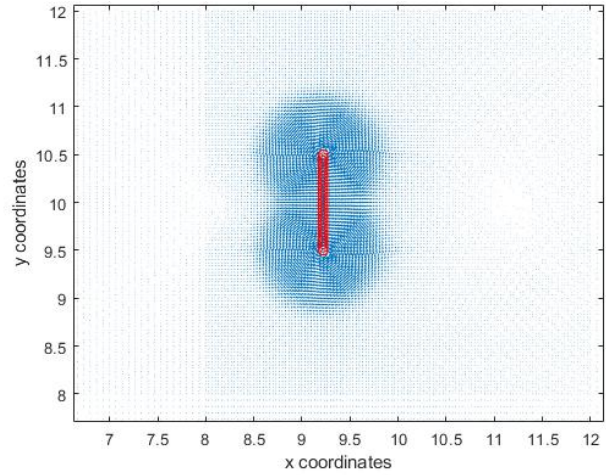

**Figure 7-4** Velocity field for $l_a \times l_b = 20B \times 20B$

## 7.2 Fully burred pipe

In this section, the thin plate studied in section 7.1 is replaced by the rigid pipe segment. The analytical solution of the non-dimensional resistance force between the fully-buried pipe and the soil is reported in (Martin and Randolph, 2006):

$$N_{pipe} = \frac{P}{Ds_u} = 11.94 \tag{7-2}$$

where $D = 1m$ is the diameter of pipe and $s_u$ is the Tresca shear stress strength.

The solution domain is $l_a \times l_b = 10D \times 10D$. The pipe cross-section is assigned with a velocity of $v_d = 0.01D$ m/s. The time increment is 0.001s. The linear Lagrangian interpolation enriched by the bubble function is used as the interpolation function for the velocity and pressure fields. There are 355 divisions in the x and y directions. The total element number is 126025. The total node number (bubble function degree-of-freedom is condensed at element level) is 126736.

The velocity field of all the four boundaries are fixed. At $x = l_a$, the pressure is also fixed as zero to provide a reference pressure.

The plot of $N_{pipe}$ computed by the proposed VMS-FEM is shown in Figure 7-5. It is observed that the numerical solution is also wave-like. The average value is very close to 12, which is very matching with the analytical solution (7-2). The velocity field is shown in Figure 7-6. The maximum shear strain rate is shown in Figure 7-7. In Figure 7-8, the maximum shear strain is shown. A clear shear band around the pipe is observed.

The pipe diameter is changed to $D = 0.5\text{m}$. The dragging velocity is still $v_p = 0.01D$ m/s. The corresponding results are shown in Figure 7-9 to Figure 7-12. From the results, it is observed that the curve of the pipe-soil resistance force verses time step also oscillate around 12, which is consistent with the analytical solution in (7-2). Some noisy points are observed. In Figure 7-9, the vibration magnitude of the pipe-soil resistance force is larger than that in Figure 7-5 because the velocity-to-diameter ratio is higher. There is more dynamical effects in the simulation of the Figure 7-11. The shear band is also clearly observed in 7-12. The shear band represents the maximum shear strain rate. It is very important for a successful numerical method to simulate this shear band, which is widely observed in many discontinuous problems (e.g. almost rigid plastic soil).

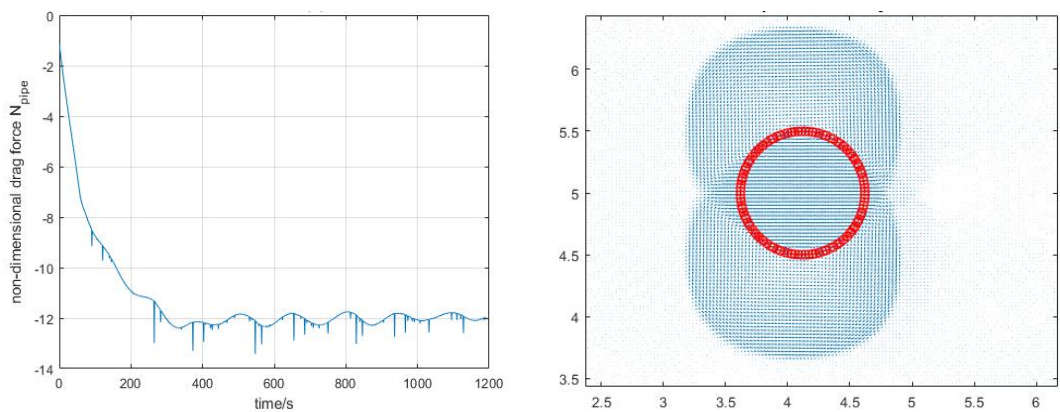

**Figure 7-5** Relationship between $N_{pipe}$ and time step; **Figure 7-6** Velocity field

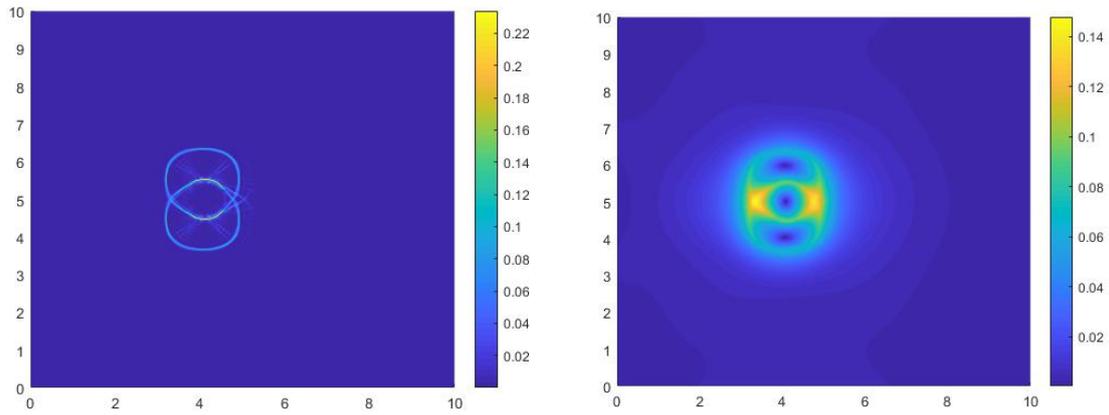

Figure 7-7 Maximum shear strain rate; Figure 7-8 Maximum shear strain

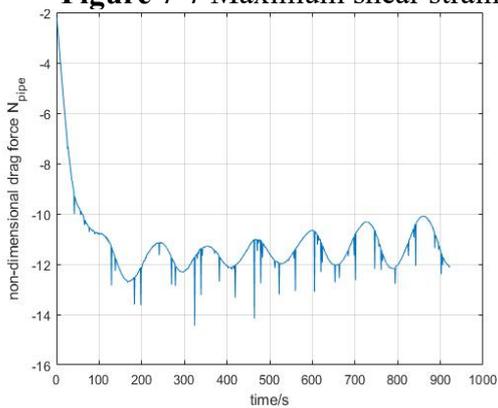 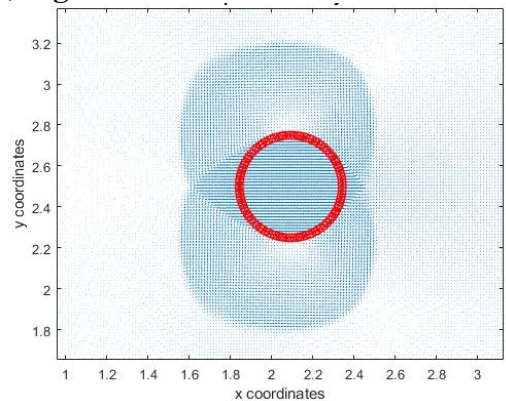

Figure 7-9 Relationship between $N_{pipe}$ and time step; Figure 7-10 Velocity field

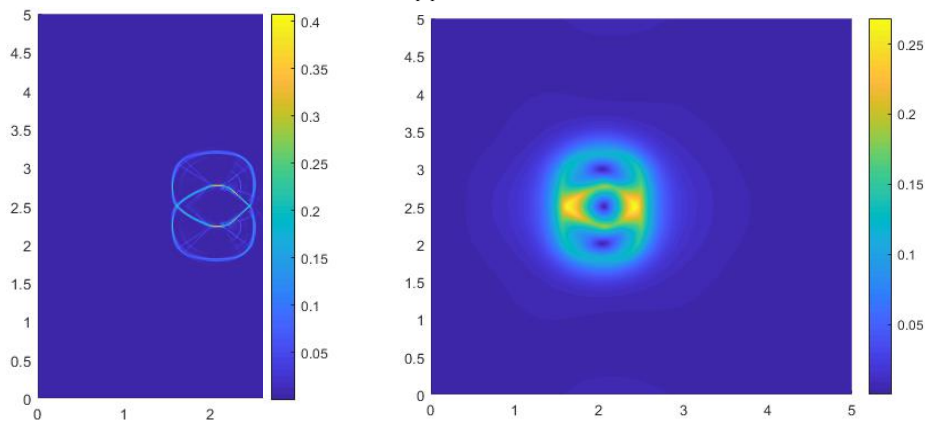

Figure 7-11 The maximum shear strain rate; Figure 7-12 The maximum shear strain

## 7.3 Pipe-soil interaction

### 7.3.1 Vertical penetration of pipe into the ideal soil

In this section, the pipe-soil vertical penetration problem is studied. The problem numerical model is shown in Figure 7-13. The water is used as the background material (blue region). The soil is the yellow region. The uniform soil is studied. No gravity is applied and the soil strength is uniform in the whole solution domain. The softening and rate dependency are also not considered in this section.

The problem's input physics parameters are:

Diameter of pipe is: $D = 1m$

Thickness of pipe is: $t_p = 0.01m$

Undrained soil strength is: $s_u = 1100 Pa$

Soil Poison ratio (volume conservative) is: $\mu_{soil} = 0.495$

Soil Von Mises yielding stress is: $\sigma_0 = 1.732 s_u$

Soil Young's Modulus is: $E_{soil} = 400 s_u$

The non-dimensional vertical resistance force is plotted against the vertical non-dimensional penetration in Figure 7-14. Noises are still observed. The ABAQUS CEL analysis results are shown in Figure 7-15. Through the comparison, it is concluded that the current non-dimensional resistance force is very accurate when comparing to the ABAQUS CEL solution.

The Tresca stress distribution is show in Figure 7-16. It is observed that the soil near the pipe is fully yielded. The Tresca stress is the same with $s_u$. The soil berm is shown in Figure 7-17. A localised soil berm is observed near the pipe.

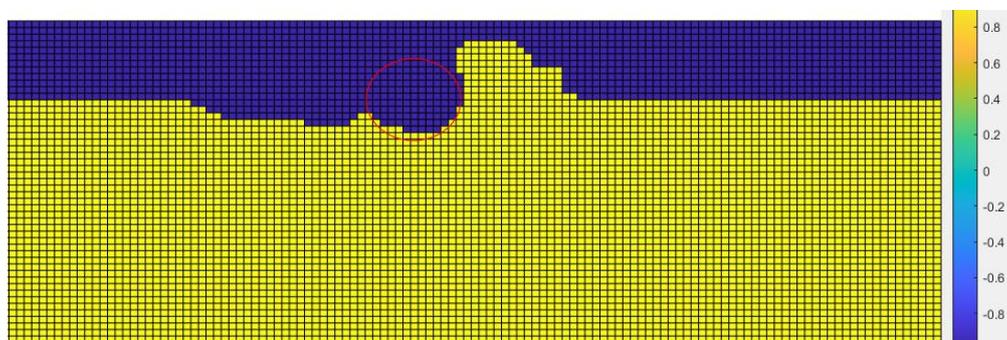

**Figure 7-13** The pipe-soil interaction problem model (Yellow: soil; Blue: water)

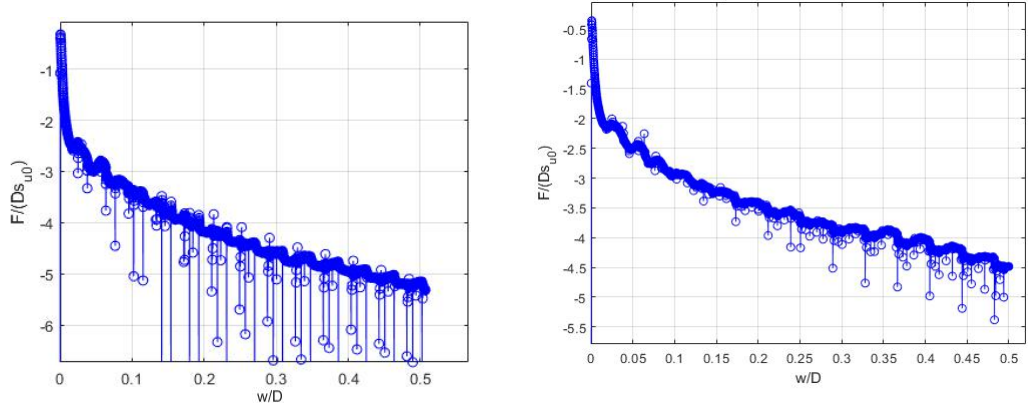

**Figure 7-14 (a, b)** Non-dimensional vertical resistance force versus non-dimensional vertical penetration $w/D$ (fully-rough: left; smooth: right)

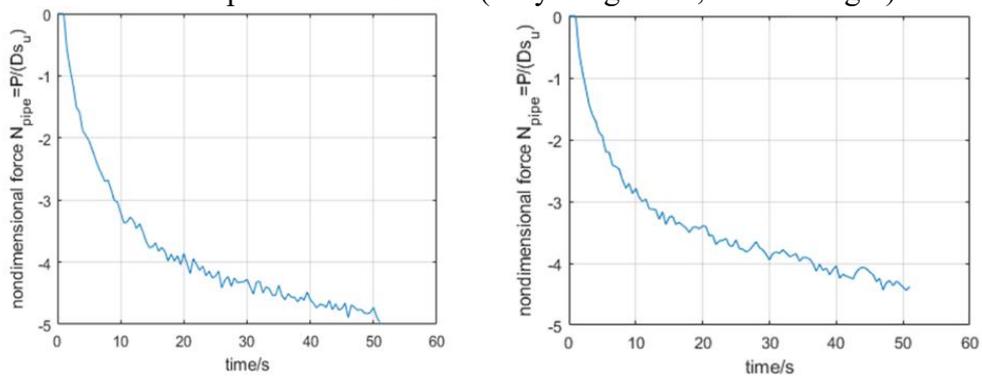

**Figure 7-15 (a, b)** Non-dimensional resistance against time (s) from the ABAQUS CEL (fully rough: left; smooth: right) $v_p = -0.01 D/s$

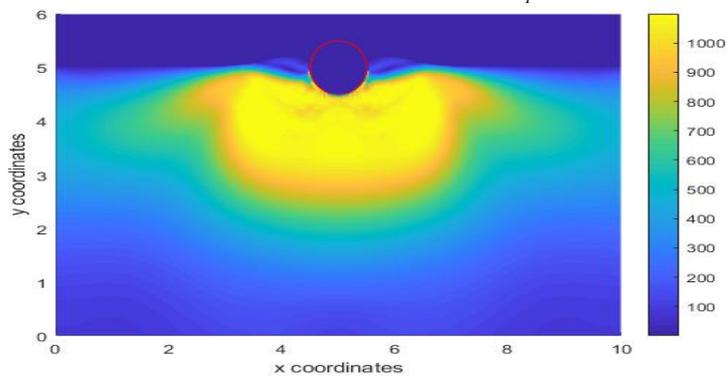

**Figure 7-16 (a)** Distribution of Tresca stress (fully-rough)

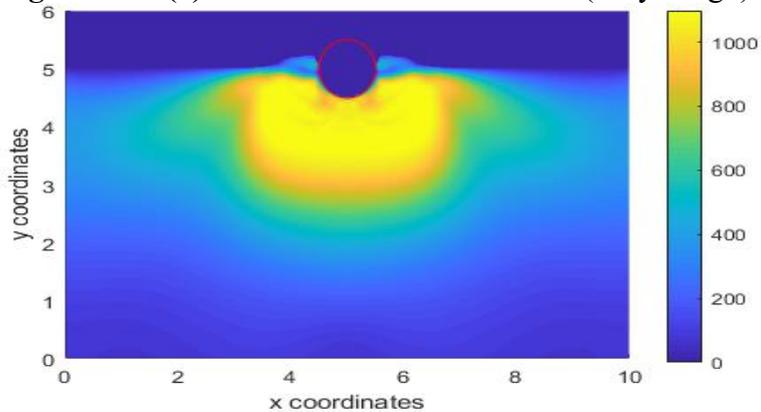

**Figure 7-16 (b)** Distribution of Tresca stress (smooth)

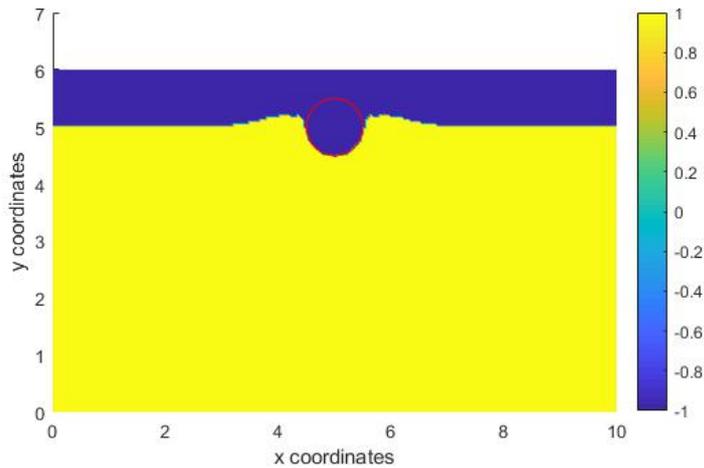
**Figure 7-17 (a)** Distribution of materials (fully-rough)

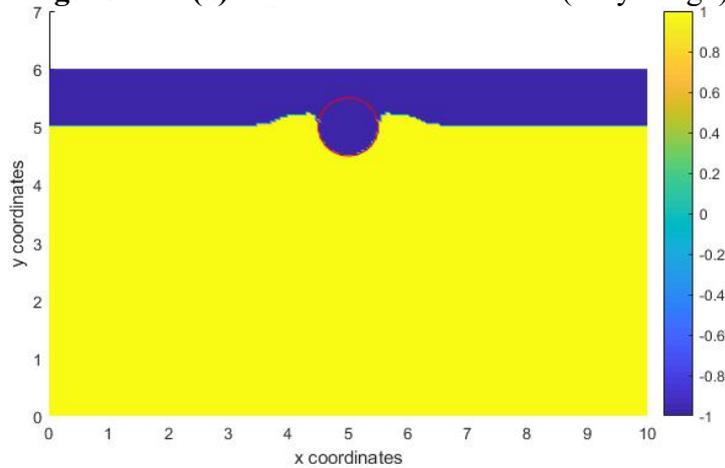
**Figure 7-17 (b)** Distribution of materials (smooth)

### 7.3.2 Horizontal movement of pipe interacting the ideal soil

In this section, the lateral movement is initiated following the vertical penetration as appeared in section 7.3.1. During the lateral movement, the vertical position of the pipe segment is fixed. The lateral velocity of the pipe segment is prescribed. The resistance force in the horizontal and vertical directions are calculated.

The loading process is detailed. A vertical velocity is assigned on the pipe to simulate the pipe-soil vertical penetration. After the total vertical penetration reaches a designed level (e.g. 0.5D), a new horizontal velocity is assigned on pipe segment to initiate the lateral movement, during which the vertical position of pipe is fixed. The non-dimensional pipe-soil resistance force during this process is studied. Both fully-rough and smooth pipe-soil friction cases are considered. Gravity and non-gravity cases are both considered. Both of the uniform soil strength and the depth-dependent soil strength are included.

The results of the pipe-soil lateral resistance force are firstly compared with those from the ABAQUS CEL in Figure 7-18. In this simulation, an initial vertical penetration of $0.5D$ is applied at first. The vertical velocity of pipe is $-0.015D$/s. Then, the vertical position of the pipe is fixed and a horizontal velocity of $0.015D$/s is applied on the pipe. The non-dimensional lateral pipe-soil resistance force $F/Ds_u$ is accurate in comparison with ABAQUS CEL method when the lateral movement is stabilized ($t > 60s$). In the VMS-FEM simulation, there is a singularity point when the pipe turns its velocity direction from vertical to lateral. The Tresca stress distribution is shown in Figure 7-19. The soil near the pipe is fully yielded. The soil berm is shown in Figure 7-20. Because no gravity is applied, the soil flow upward in Figure 7-20.

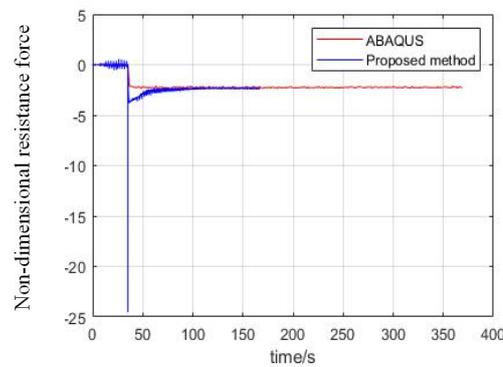

**Figure 7-18** The non-dimensional resistance force versus time (fully rough friction)

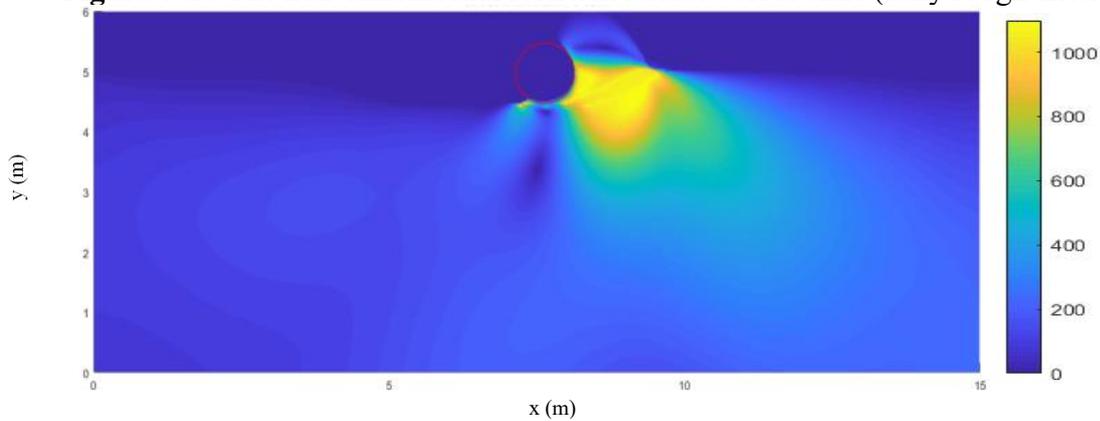

**Figure 7-19** The Tresca stress distribution

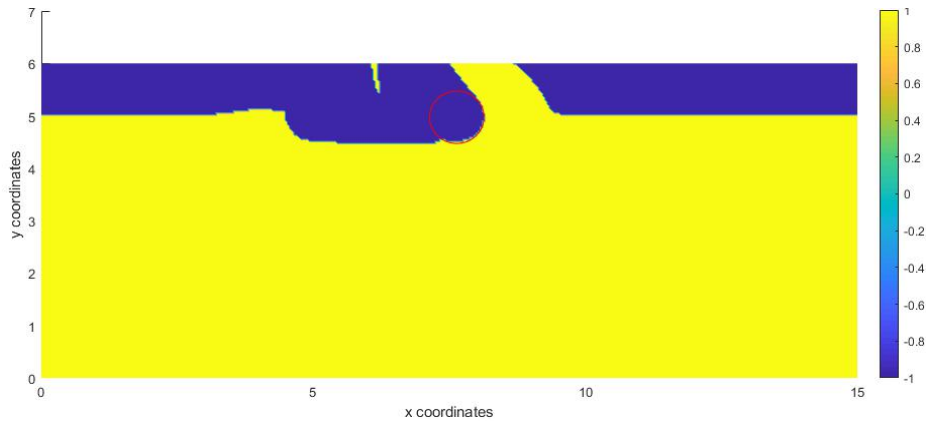

**Figure 7-20** Distribution of materials during lateral movement

### 7.3.3 Vertical penetration of pipe into the realistic soil

The pipe-soil penetration problem with existence of soil gravity, soil strain-rate-dependent hardening effect and the soil plastic-strain-dependent softening effect is analysed in this section. The problem model is the same with that in section 7.3.1. The only differences are that the gravity of soil is applied; the strain-rate-dependent hardening constitutive model is considered with the plastic-strain-dependent softening effect, see (6-3). The input physics parameters are given in Table 7-1 (Kong, 2015). The non-dimensional resistance force results for the ideal soil case are compared in Figure 7-21. The soil berm configuration and the Tresca stress distribution are shown in Figure 7-22 and Figure 7-23.

The convergence study on the mesh density are presented in Figure 7-24. The strain rate dependent soil problems are studied, see Figure 7-21 b for the three cases of $\mu = 0, 0.05, 0.1$. Both of the fully-rough and smooth pipe-soil friction models are considered. By changing the element size from 0.0533m to 0.032m, the results are almost unchanged.

**Table 7-1** The physics parameters of the problem (Kong, 2015)

| Parameter | Value |
|---|---|
| Pipe diameter, $D$: m | 0.8 |
| Shear strength of soil at mudline, $s_{um}$: kPa | 2.3 |
| Shear strength gradient, $k$: kPa/m | 3.6 |
| Submerged unit weight of soil, $\gamma'$: kN/m$^3$ | 6.5 |
| Initial pipe embedment, $w_{ini}$ | $0.45D$ |
| Lateral pipe displacement, $u_{fin}$ | $3D$ |
| Pipe weight during lateral movement, $W$: kN/m | 3.3 |
| Vertical penetration rate, $v_p$ | $0.015D$/s |
| Lateral sweep rate, $v_p$ | $0.050D$/s |

| Parameter | Value |
|---|---|
| Rate of strength increase per decade, $\mu$ | 0.10 |
| Reference shear strain rate, $\dot{\gamma}_{ref}$: s$^{-1}$ | $3\times10^{-6}$s$^{-1}$ |
| Sensitivity of soil, $S_t$ | 3.2 |
| Ductility parameter of soil, $\xi_{95}$ | 10 |

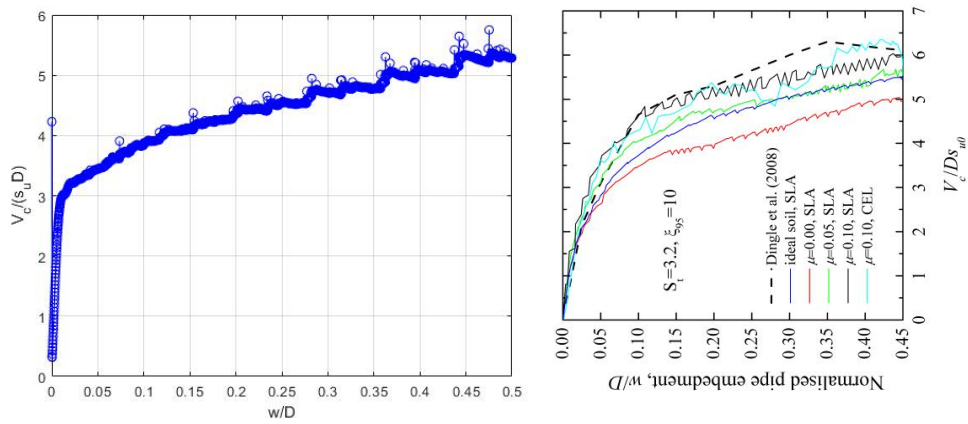

**Figure 7-21** The non-dimensional resistance force for the ideal soil case (**a**: VMS-FEM; **b**: the blue curve from the SLA by Kong (2015))

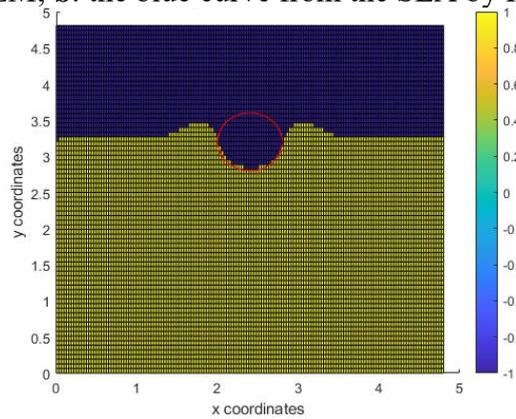

**Figure 7-22** Distribution of materials of water and soil

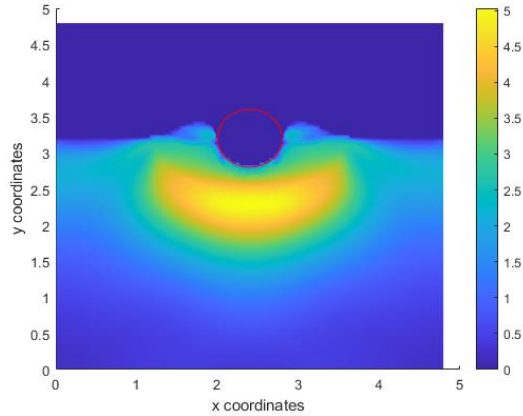

**Figure 7-23** The distribution of Tresca stress

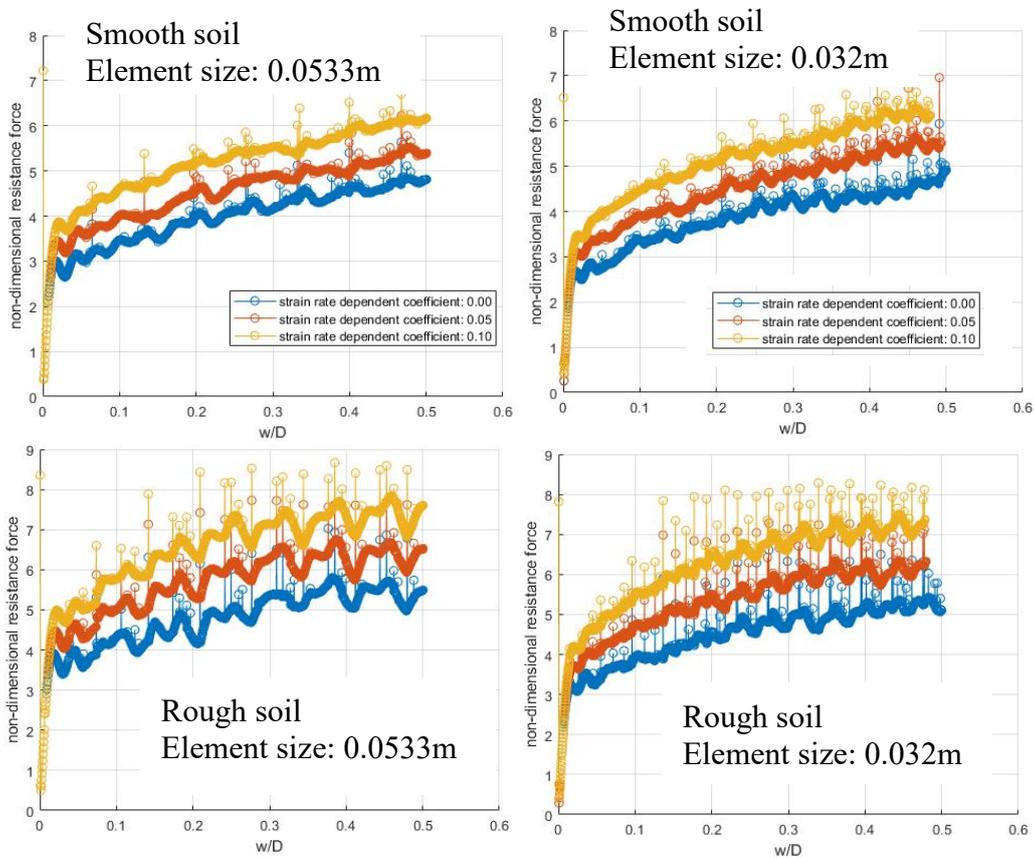

**Figure 7-24** Convergence study on mesh densities

### 7.3.4 Horizontal movement of pipe interacting the realistic soil

The lateral movement of the pipe interacting with the soil foundation under the gravity is studied. The outer diameter of pipe is $D = 0.8m$. An initial vertical penetration of 0.1D, 0.2D, 0.3D, 0.4D, 0.5D with the vertical velocity -0.01 D/s is applied at first. Then, after 10s relaxation, a horizontal velocity 0.01 D/s is assigned to the pipe to initiate the lateral movement. The solution domain must be large enough in the horizontal direction. In this simulation, we define 9D width for the solution domain. The soil property and the pipe property are the same with that in section 7.3.3.

For the case without the plastic-strain-dependent softening effect and the strain-rate-dependent hardening effect, the non-dimensional pipe-soil resistance force is plotted in Figure 7-25. The soil berm is shown in Figure 7-26. In Figure 7-27, the distribution of the Tresca stress is presented. The non-dimensional soil strength gradient $k = |y_0 - y_0^{free}| k_{su} / s_u^0$ is shown in Figure 7-28, where $y_0$ is the initial vertical coordinate of the soil material points, $y_0^{free} = 3.2m$ is the initial vertical coordinate of the soil free surface, $s_u^0 = 2.3kPa$ is the free surface soil strength, $k_{su} = 3.6kPa/m$ is the initial soil gradient along depth, see Table 7-1.

It is observed that the vertical and horizontal resistance forces increase with the vertical penetration. The accumulated soil berm also becomes larger with the vertical penetration. The use of the bubble function in the velocity interpolation function makes the calculated pipe-soil resistance force results smaller than those from the simulations without using the bubble function. However, the difference is not very large.

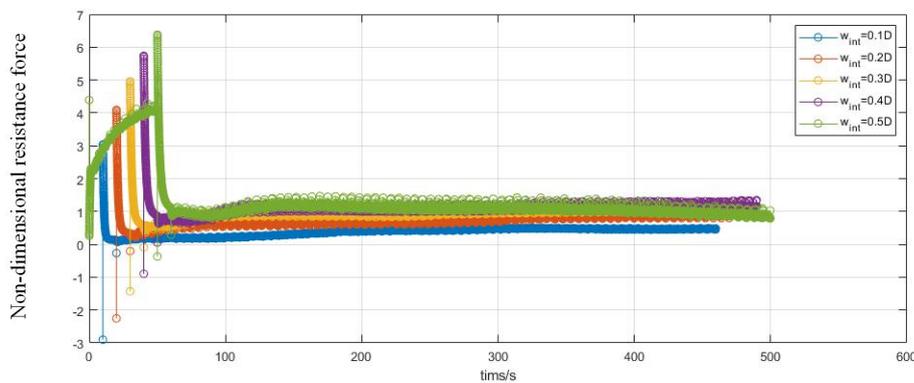

**Figure 7-25 (a)** The non-dimensional vertical resistance force (bubble function)

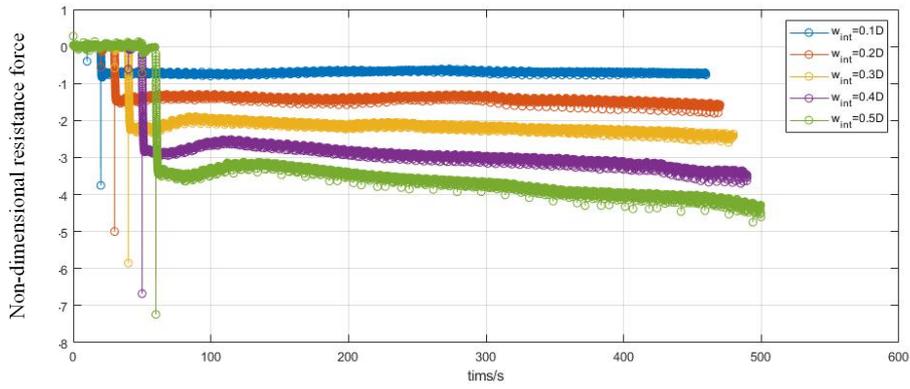

**Figure 7-25 (b)** The non-dimensional horizontal resistance force (bubble function)

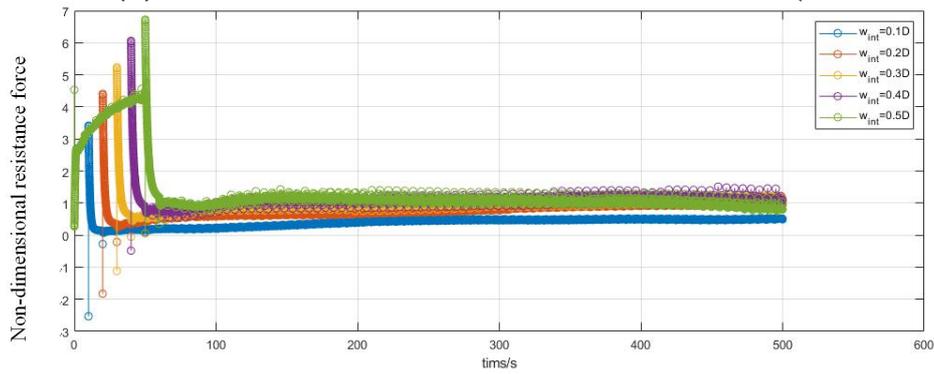

**Figure 7-25 (c)** The non-dimensional vertical resistance force (no bubble function)

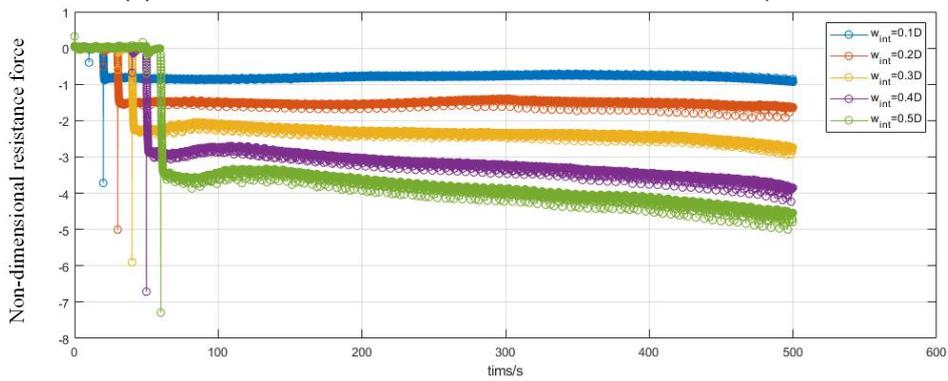

**Figure 7-25 (d)** The non-dimensional horizontal resistance force (no bubble function)

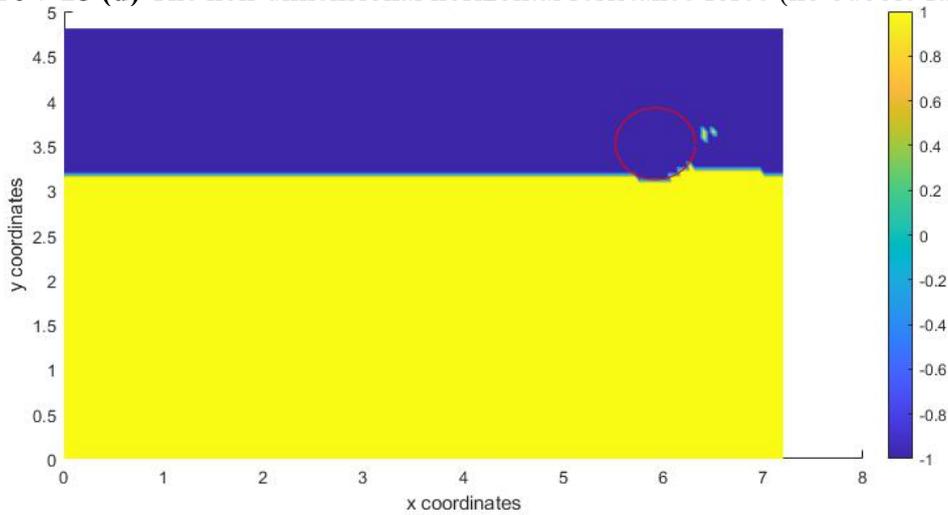

**Figure 7-26 (a)** Distribution of materials (0.1D initial vertical penetration)

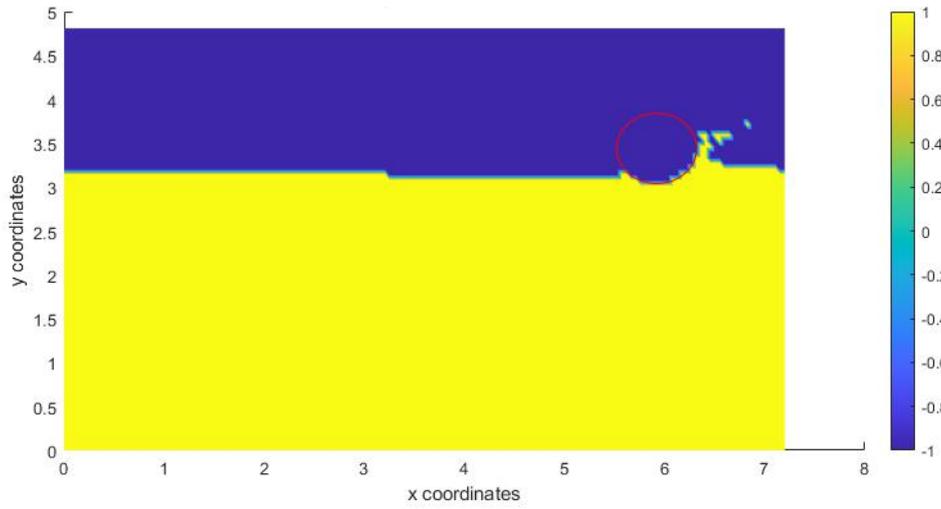

**Figure 7-26 (b)** Distribution of materials (0.2D initial vertical penetration)

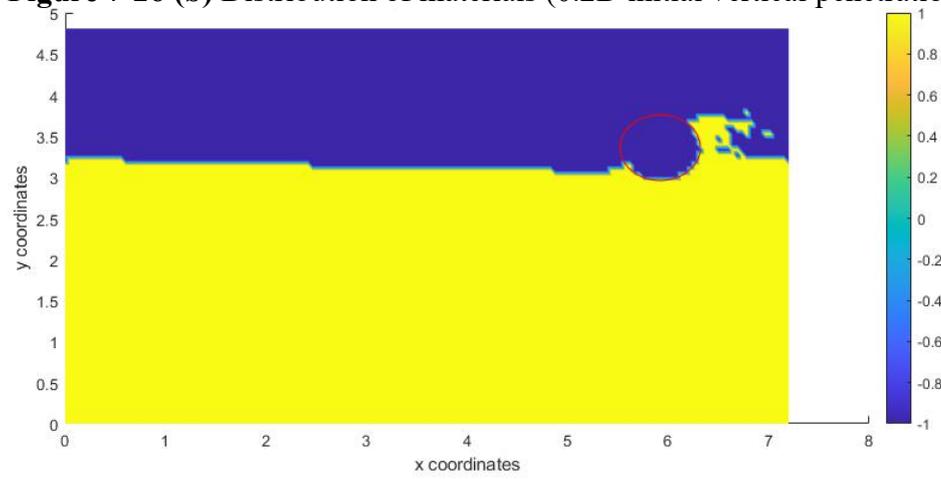

**Figure 7-26 (c)** Distribution of materials (0.3D initial vertical penetration)

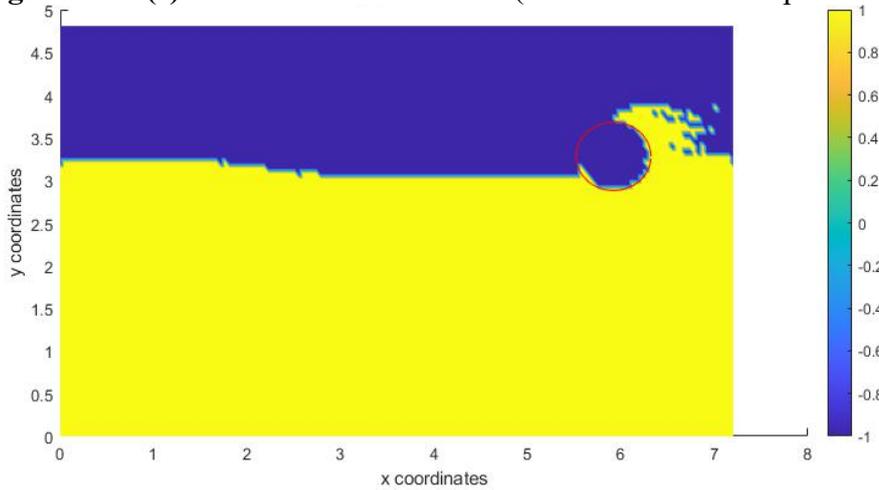

**Figure 7-26 (d)** Distribution of materials (0.4D initial vertical penetration)

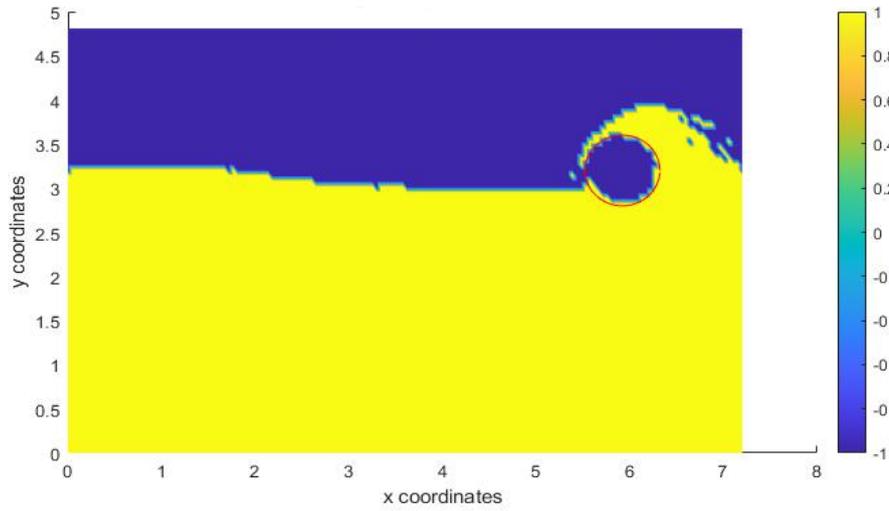

**Figure 7-26 (e)** Distribution of materials (0.5D initial vertical penetration)

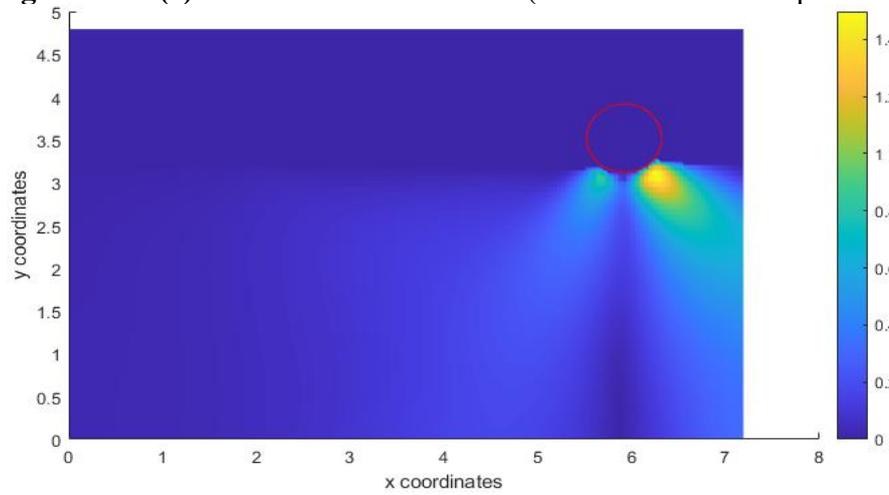

**Figure 7-27 (a)** Distribution of Tresca stress (kPa) (0.1D initial vertical penetration)

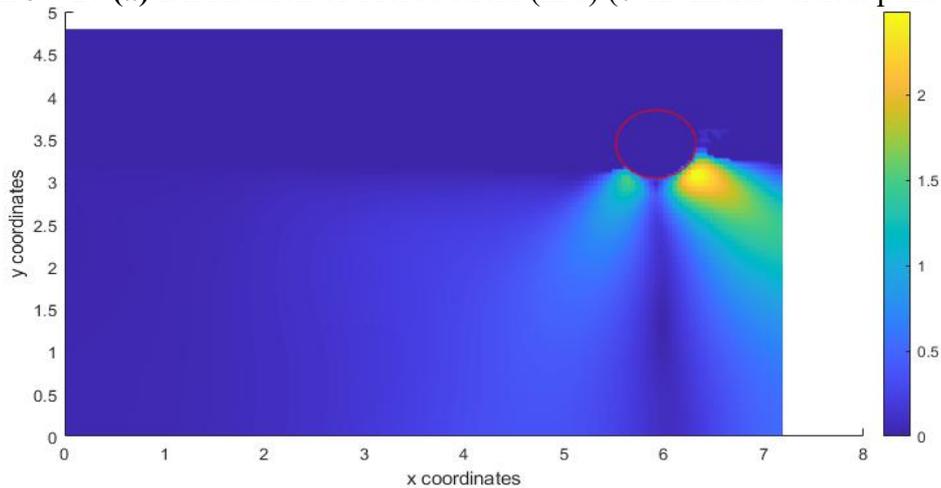

**Figure 7-27 (b)** Distribution of Tresca stress (kPa) (0.2D initial vertical penetration)

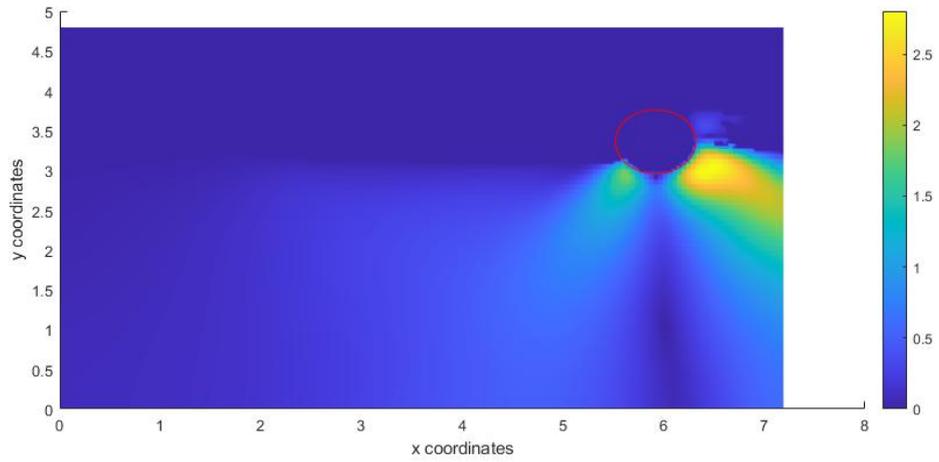

**Figure 7-27 (c)** Distribution of Tresca stress (kPa) (0.3D initial vertical penetration)

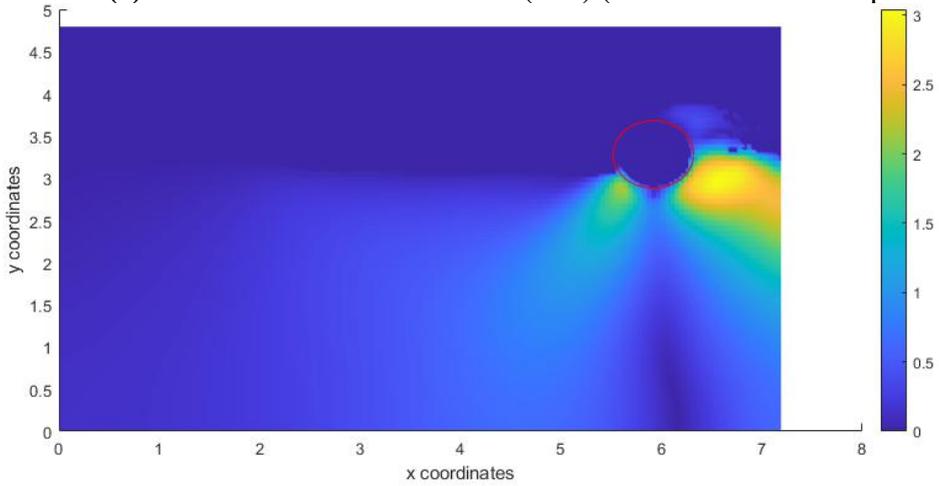

**Figure 7-27 (d)** Distribution of Tresca stress (kPa) (0.4D initial vertical penetration)

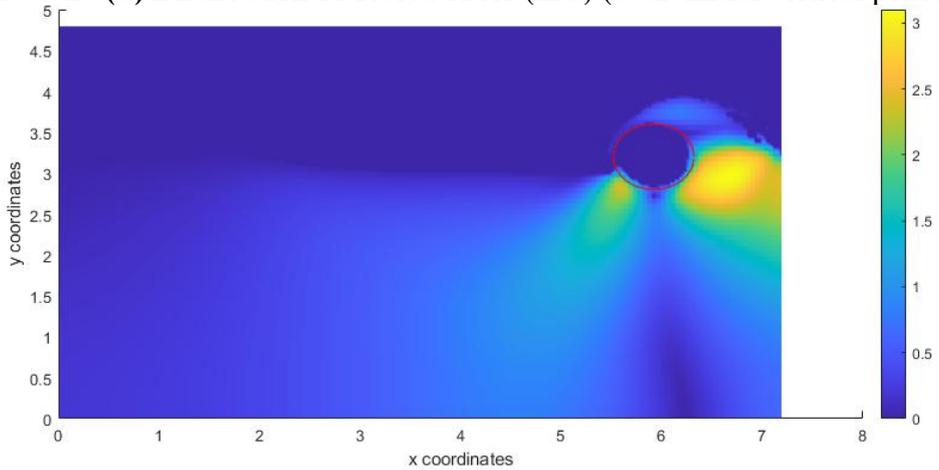

**Figure 7-27 (e)** Distribution of Tresca stress (kPa) (0.5D initial vertical penetration)

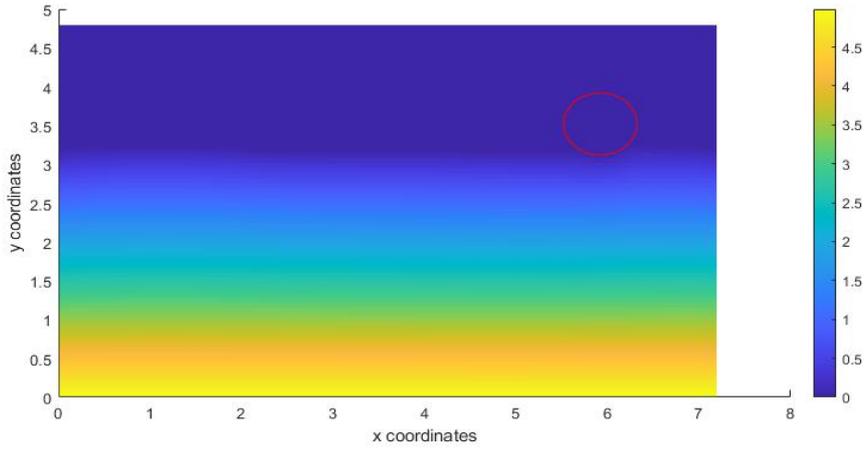

**Figure 7-28 (a)** Distribution of non-dimensional strength gradient
$k = \left| y_0 - y_0^{free} \right| k_{su} / s_u^0$ (0.1D initial penetration)

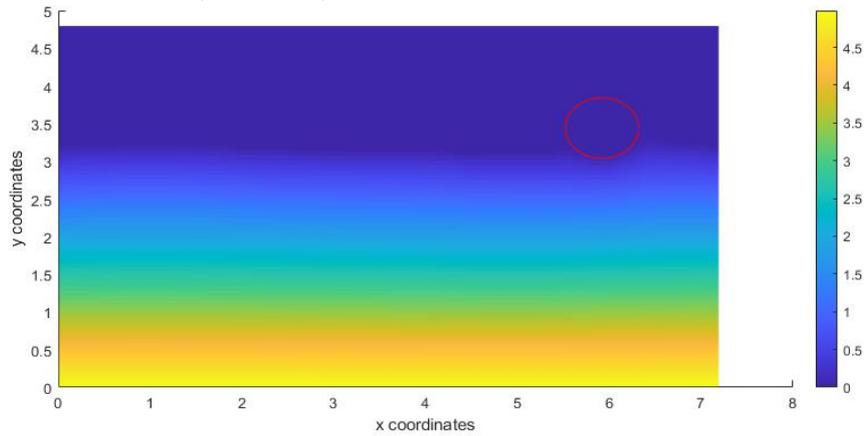

**Figure 7-28 (b)** Distribution of non-dimensional strength gradient
$k = \left| y_0 - y_0^{free} \right| k_{su} / s_u^0$ (0.2D initial penetration)

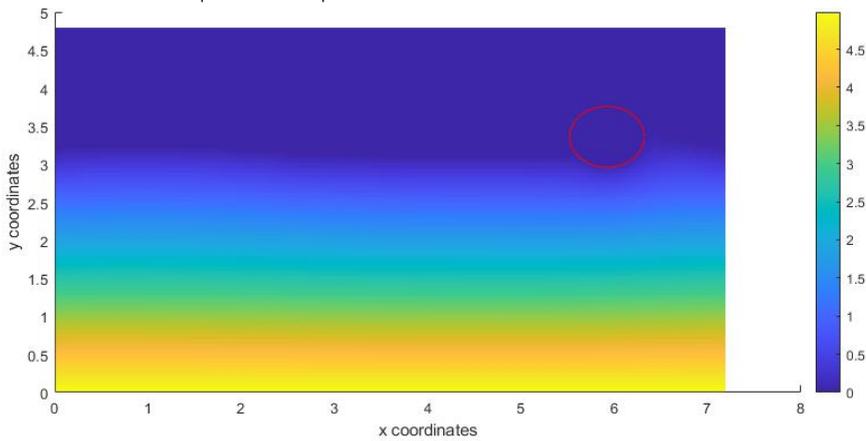

**Figure 7-28 (c)** Distribution of non-dimensional strength gradient
$k = \left| y_0 - y_0^{free} \right| k_{su} / s_u^0$ (0.3D initial penetration)

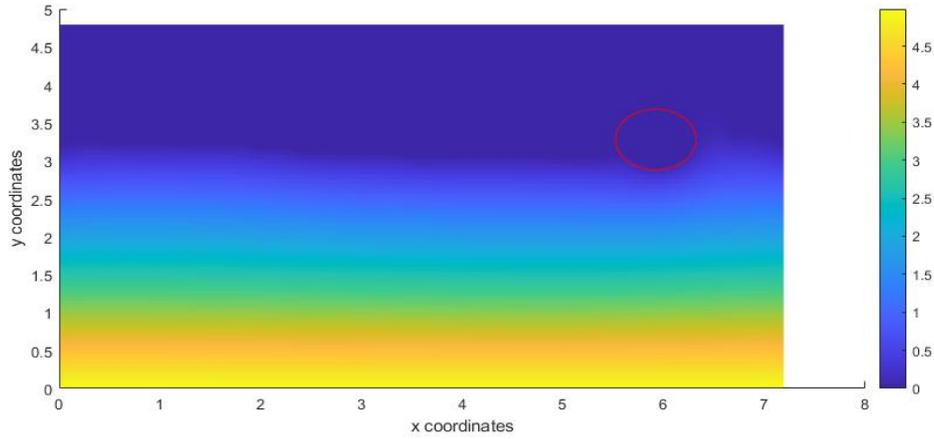

**Figure 7-28 (d)** Distribution of non-dimensional strength gradient
$k = \left| y_0 - y_0^{free} \right| k_{su} / s_u^0$ (0.4D initial penetration)

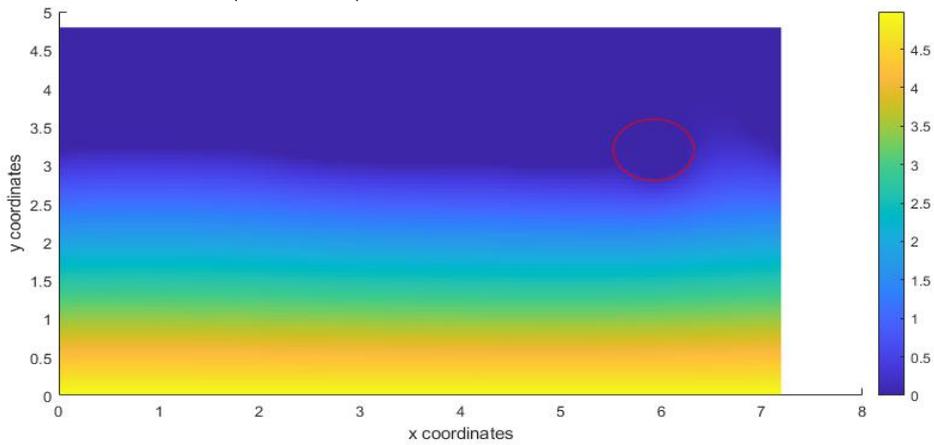

**Figure 7-28 (e)** Distribution of non-dimensional strength gradient
$k = \left| y_0 - y_0^{free} \right| k_{su} / s_u^0$ (0.5D initial penetration)

The soil with the softening effect and strain rate hardening effect is studied. The operative soil strength is a function of the strain rate and the total plastic strain, see equation (6-3). The corresponding results of the pipe-soil resistance forces are shown in Figure 7-29. The distribution of the current non-dimensional soil strength gradient $wk_{su} / s_u^0$ are shown in Figure 7-30. The soil berm is shown in Figure 7-31. It is observed that the pipe-soil resistance force in Figure 7-29 is slightly smaller than that in Figure 7-25 because the softening effect due to the accumulated plastic strain is more dominating than the hardening effect due to the strain rate.

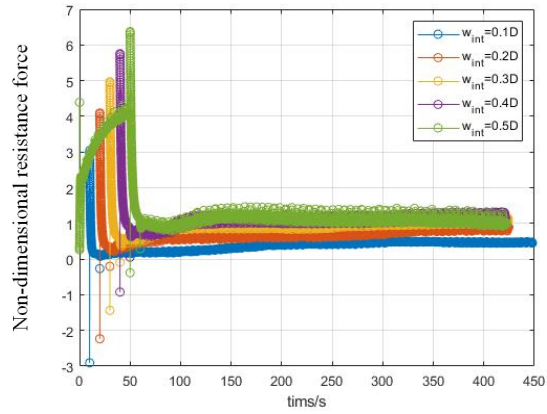

**Figure 7-29 (a)** The non-dimensional vertical resistance force (bubble function)

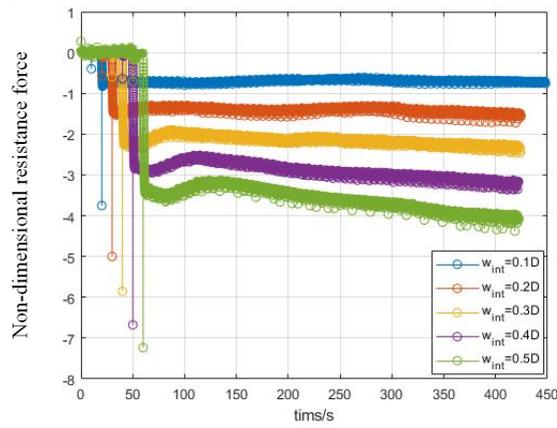

**Figure 7-29 (b)** The non-dimensional horizontal resistance force (bubble function)

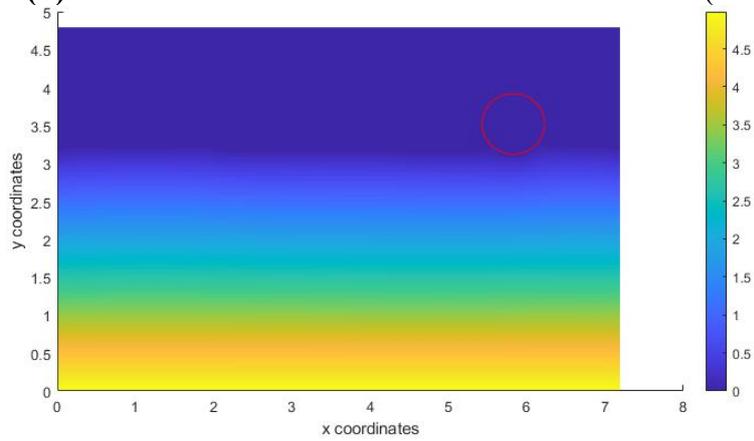

**Figure 7-30 (a)** Distribution of non-dimensional strength gradient
$k = \left| y_0 - y_0^{free} \right| k_{su} / s_u^0$ (0.1D initial penetration)

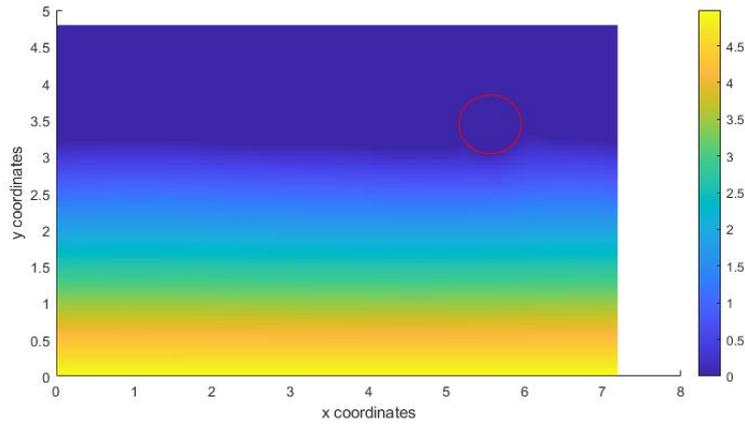

**Figure 7-30 (b)** Distribution of non-dimensional strength gradient
$k = |y_0 - y_0^{free}| k_{su} / s_u^0$ (0.2D initial penetration)

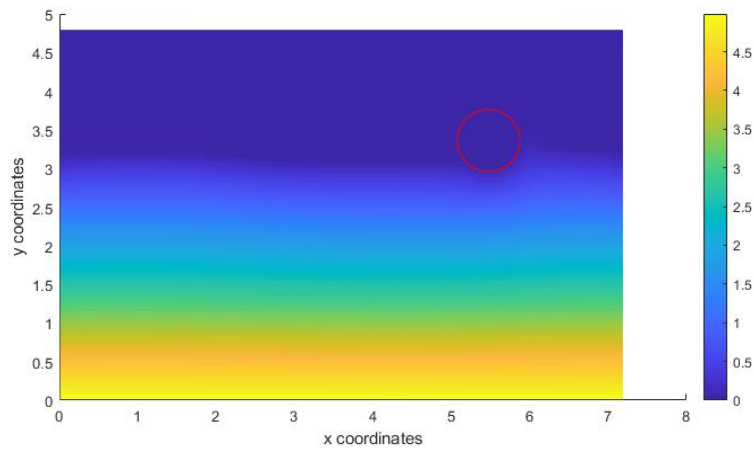

**Figure 7-30 (c)** Distribution of non-dimensional strength gradient
$k = |y_0 - y_0^{free}| k_{su} / s_u^0$ (0.3D initial penetration)

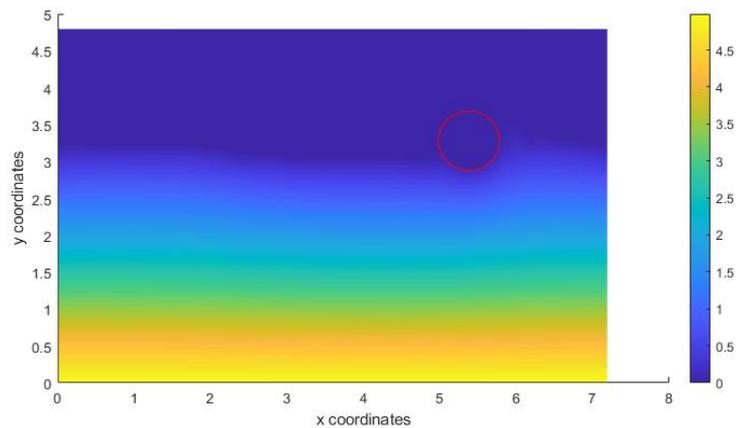

**Figure 7-30 (d)** Distribution of non-dimensional strength gradient
$k = |y_0 - y_0^{free}| k_{su} / s_u^0$ (0.4D initial penetration)

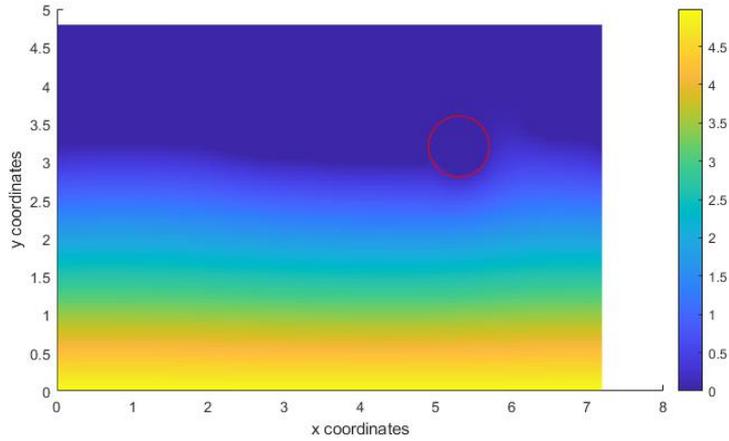

**Figure 7-30 (e)** Distribution of non-dimensional strength gradient
$$k = \left|y_0 - y_0^{free}\right| k_{su} / s_u^0 \quad \text{(0.5D initial penetration)}$$

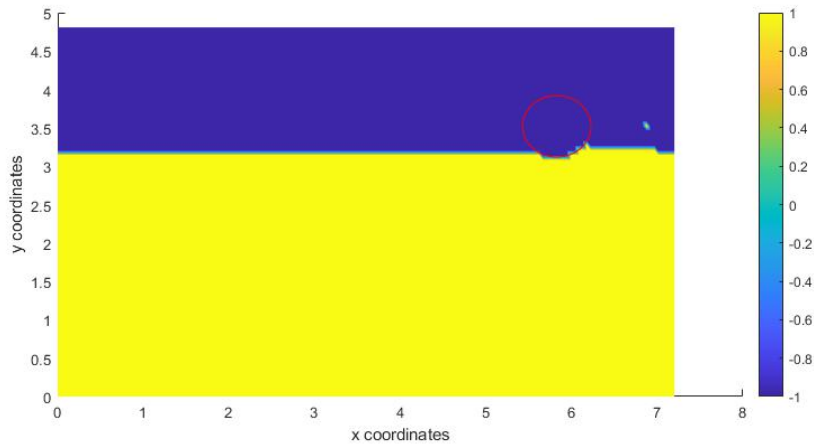

**Figure 7-31 (a)** Distribution of materials (0.1D initial vertical penetration)

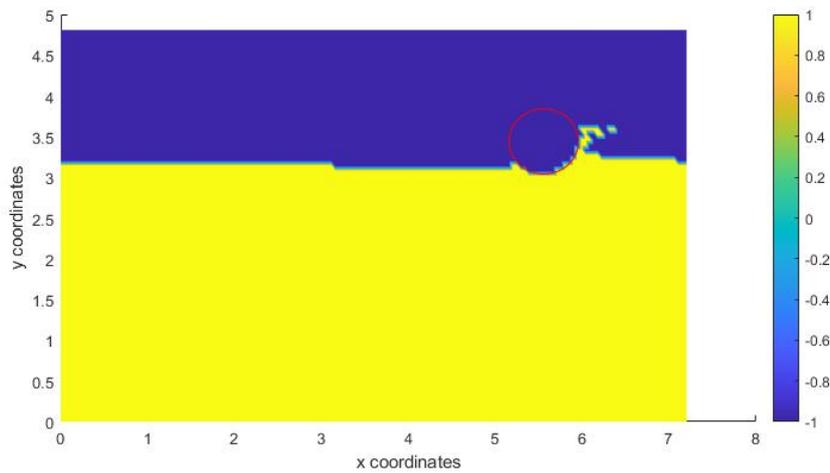

**Figure 7-31 (b)** Distribution of materials (0.2D initial vertical penetration)

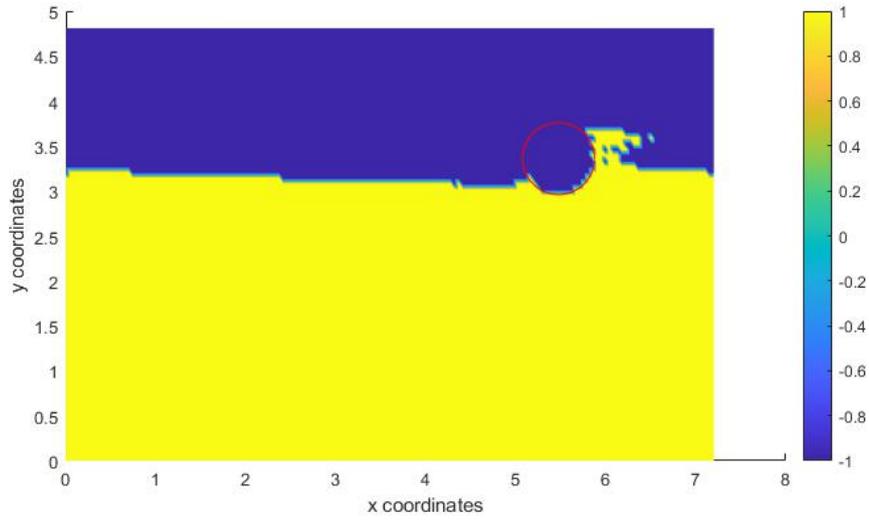

**Figure 7-31 (c)** Distribution of materials (0.3D initial vertical penetration)

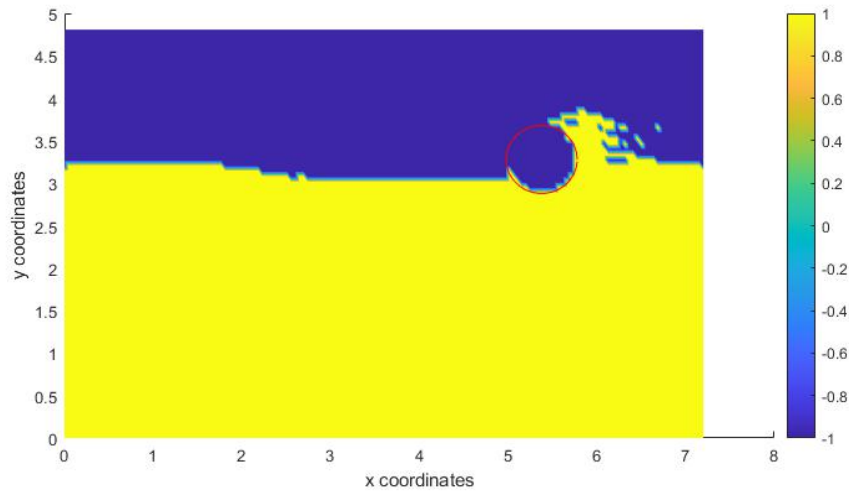

**Figure 7-31 (d)** Distribution of materials (0.4D initial vertical penetration)

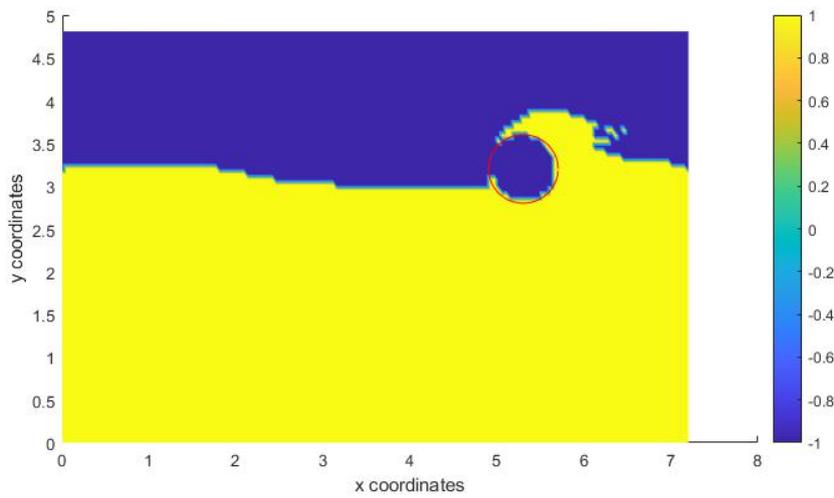

**Figure 7-31 (e)** Distribution of materials (0.5D initial vertical penetration)

## 8. Conclusion

The proposed VMS-FEM has been successfully applied to study the history dependent material in the Eulerian framework. The pipe-soil interaction analysis is focused. The accuracy and the stability are both good.